# Poorly Separated Infinite Normal Products

## N. Noble[*]

**Abstract**

From the literature: a product of compact normal spaces is normal; the product of a countably infinite collection of non-trivial spaces is normal if and only if it is countably paracompact and each of its finite sub-products is normal; if all powers of a space X are normal then X is compact – provided in each case that the spaces involved are $T_1$. Here I examine the situation for infinite products not required to be $T_1$ (or regular), extending or generalizing each of these facts. In addition, I prove some related results, give a number of examples, explore some alternative proofs, and close with some speculation regarding potential applications of these findings to category theory and lattice theory. I plan to consider the normality of poorly separated finite products in a companion paper currently in preparation.

**Summary of principal results**

Calling a space ***vacuously normal*** if it does not contain a pair of nonempty disjoint closed subsets, and using the notation $\Pi_{\alpha \in S} X_\alpha = X_S$, I show:

- A product of vacuously normal spaces is vacuously normal.
- A product of compact normal spaces is compact normal.
- All powers of X are normal if and only if X is vacuously normal or compact normal.
- If X is vacuously normal, Y normal, and the projection from X x Y to Y is closed, then XxY is normal; the converse holds for Y a $T_1$ space, but not in general.
- A countable product with no factor vacuously normal is normal if and only if it is countably paracompact and each of its finite sub-products is normal.
- If $X_A$ is normal then the index set A can be expressed as the union of disjoint subsets B, C, D, and E (some possibly empty) where:
  - $X_B$ is vacuously normal;
  - $X_C$ is compact;
  - $X_D$ is $[\omega, \kappa]$-compact, where $\kappa = |A \setminus B|$; and
  - $X_E$ is $\kappa$-paracompact with E countable.

[*] P.O. Box 1812, Los Gatos, CA 95031; norm@metavideo,com

# 0  Introduction

To be clear, a space is **normal** if disjoint closed sets are contained in disjoint open sets and **normal $T_1$** if in addition it is regular or $T_1$ (and therefore completely regular Hausdorff, *i.e.*, Tychonoff). The term $T_4$ is often used for one or the other of normal / normal $T_1$, but with no consistency among authors.

With no additional separation, normal spaces still enjoy what are arguably their most distinguishing properties - those captured by Urysohn's Lemma and the Tietz Extension Theorem. But when it comes to discussing products, additional separation is almost always assumed (if not already present in the author's definition of a "space"). This is true of Kelley [1955], Willard [1970], Engelking [1989], as well as many other texts and monographs dealing with general topology ranging from Hausdorff [1914] to Goubault-Larecq [2013]; the treatise on normality by Alo & Shapiro [1974] (which, in fact, does not treat products); surveys including Przymusinski [1984], Atsuji [1989], and Hoshina [1989]; updates by Gruenhge [1992], [2002]; and numerous proceeding papers or journal articles. There are a few exceptions: in particular, treatments involving locally finite families (*e.g.*, papers by Dowker, Morita, and Katuta ).

While historically point set topologists have favored Hausdorff or stronger separation, non-Hausdorff topologies frequently arise in association with other structures. These include the Zariski topologies of algebraic geometry, the digital topologies of computer graphics, and the various topologies, notably the upper $\subseteq$ Scott $\subseteq$ Lawson $\subseteq$ Alexandroff topologies associated with partial orders and lattices.[1]

There is a rich literature concerning non-Hausdorff spaces and their relationship to orders, surveyed and synthesized in Gierz, Hofmann, Keimel, Lawson, Mislov, & Scott [1980] and [2003], and Goubault-Larecq [2013]. However, to date I have found little of significance in that literature relating to normality of non-$T_1$ products. (Works which consider interval topologies - LOTs, PO-spaces, GO-spaces - have considered normality from the beginning, Nachbin [1950], but such spaces are Hausdorff.)

# 1  Terminology, notation, and preliminary results

Here I establish some conventions, mostly standard; identify some spaces which will prove useful; provide some examples which may be informative; and prove some preliminary results in preparation for the later sections. For $X=\Pi_{\alpha \in A}X_\alpha$ and $B \subset A$, write $X = X_A$, $X_B = \Pi_{\alpha \in B}X_\alpha$, and $X_{A \setminus B} = \Pi_{\alpha \in A \setminus B}X_\alpha$, so for example, $X = X_B \times X_{A \setminus B}$. If an index set I is empty, $X_I$ is understood to be the space with a single point.

---

[1] Associated with any topology is its order of specialization: $x \leq y \Leftrightarrow y \in \bar{x}$ (*i.e.*, iff each open set containing y contains x). The upper topology is the coarsest one producing a given order as its order of specialization; it takes the closures of points as a subbase for the closed sets. The Alexandroff topology is the largest yielding the given order, taking as closed all unions of point closures (so the intersection of any family of open sets is open and the space is "Alexandroff discrete").

## 1.1 Numbers and elementary spaces

- **0,1**    when used as logic, false / true respectively (thus below "$s_{\alpha\beta} = (\alpha > \beta)$" is used to describe the zero / one coordinate of a point in an Alexandroff or Cantor cube.
- **II**    $\{0,1\}$ with topology $\{\varnothing, \mathbf{II}\} \cup \{\mathbf{0}\}$, the Serpinski space.
- **2**    $\{0,1\}$ with the discrete topology.
- **E**    $\{0,1,2\}$ with topology $\{\varnothing, \mathbf{E}\} \cup \{\mathbf{0}\}$, (see footnote in next subsection).
- **ℵ**    an infinite cardinal and a set of that cardinality with, when bold, the discrete topology, $\mathbf{ℵ_0}$ and $\mathbf{ℵ_1}$ being the first and second infinite cardinals, but
- **N**    the countably infinite set and often **N** instead of $\mathbf{ℵ_0}$ for the discrete space.
- **ω**    the first infinite ordinal, the cardinal $ℵ_0$, or the ordered space $[0,\omega)$.
- **$\overline{\omega}$**    (or sometimes $\mathbf{N^*} = \mathbf{N} \cup \{\infty\}$ or $\{0\} \cup \{1/n: n > 0\}$): $[0,\omega]$ with the order topology;

I make similar use of $\Omega$, $\overline{\Omega}$, for the first uncountable ordinal and, *e.g.*, $\kappa$ and $\overline{\kappa}$ for limit ordinals / cardinals of unspecified size. Note that $\overline{\kappa}$ is both the Alexandroff "one point" compactification and, for $\kappa > \omega$, the Čech-Stone compactification of $\kappa$. I will have occasion to consider $\kappa = [0,\kappa)$ and $\kappa^+ = [0,\kappa]$ with two other non-discrete topologies:

- **↓κ**    $\kappa = [0,\kappa)$ with the lower topology $\{\kappa\} \cup \{[0,\alpha): \alpha \in \kappa\}$,
- **↓$\overline{\kappa}$**    $\kappa^+ = [0,\kappa]$ with the lower topology $\{\kappa^+\} \cup \{[0,\alpha): \alpha \in \kappa\}$,
- **↑κ**    $\kappa = [0,\kappa)$ with the upper topology $\{\varnothing, \kappa\} \cup \{(\alpha.\kappa)\ \alpha \in \kappa\}$,
- **↑$\overline{\kappa}$**    $\kappa^+ = [0,\kappa]$ with the upper topology $\{\varnothing, \kappa^+\} \cup \{(\alpha.\kappa]\ \alpha \in \kappa\}$, and
- **$\overline{\kappa}(\iota)$**    $\kappa^+ = [0,\kappa]$ with the topology $\{U \subseteq \kappa^+: \kappa \notin U \text{ or } |\kappa^+ \setminus U| < \iota\}$.

## 1.2 Cubes

Cubes are useful as "test" spaces, as in "X is κ-paracompact if and only if X x $\mathbf{I}^\kappa$ is normal" or as universal spaces, as in "each Tychonoff space of weight not exceeding κ can be embedded in $\mathbf{I}^\kappa$", especially because they often appear as well understood subspaces in products of other spaces. Several cubes will appear in our discussions of normal products and for convenience I have gathered some needed information about them here, plus additional facts included purely for our viewing pleasure. All of this information can be found in Engelking [1989], along with a description of the space **J**(κ), the hedgehog of spininess κ.

As indicated above, the Greek letter κ will usually be assumed to be an infinite cardinal so each cube $\mathbf{X}(\kappa) = X^\kappa$ is of weight κ. When a specific cardinal κ is being considered the name of the cube is extended by the addition of "of weight κ". As usual, "$\mathbf{X}(\kappa)$ is universal for $\mathscr{P}$" means that each space in $\mathscr{P}$ of weight not exceeding κ can be embedded in $\mathbf{X}(\kappa)$. A space is *dyadic* if it is a continuous image of the Cantor set.

| Symbol | Space | Name | Comments |
|---|---|---|---|
| $A(\kappa)$ | $\mathbb{II}^\kappa$ | Alexandroff cube | universal for $T_0$ spaces |
| $C(\aleph_0)$ | $2^\omega$ | Cantor set | every compact metrizable space is dyadic |
| $C(\kappa)$ | $2^\kappa$ | Cantor cube | universal for extremely disconnected spaces |
| $E(\kappa)$ | $\mathbf{E}^\kappa$ | Engelking[2] cube | universal for all spaces |
| $N(\kappa)$ | $\mathbf{N}^\kappa$ | Stone[3] cube | for $\kappa > \omega$, see below for $\kappa = \omega$; $N(\kappa)$ is not normal |
| $P$ | $\mathbf{N}^\omega$ | Baire space | of weight $\aleph_0$, homeomorphic to the irrationals |
| $B(\kappa)$ | $\kappa^\omega$ | Baire space | universal for metrizable spaces X with $\mathrm{Ind}(X)=0$ |
| $I(\aleph_0)$ | $\mathbf{I}^\omega$ | Hilbert cube | universal for compact metrizable spaces; universal for separable metrizable spaces |
| $I(\kappa)$ | $\mathbf{I}^\kappa$ | Tychonoff cube | universal for Tychonoff spaces |
| $H$ | $\mathbf{R}^\omega$ | Hilbert space | of weight $\aleph_0$, useful in Quantum mechanics |
| $H(\kappa)$ | $\mathbf{J}(\kappa)^\omega$ | Hilbert space | of weight $\kappa$; universal for all metric spaces[4] |

I would be remiss if I failed to mention that $\mathbf{I}(\aleph_0) \times \mathbf{B}(\kappa)$ is universal for strongly paracompact metrizable spaces. Also, the rationals are universal for the class of countable metrizable spaces, but are not much of a cube.

Since $\mathbf{I}$ is a continuous image of $\mathbf{C}(\aleph_0)$, $\mathbf{I}(\kappa)$ is a continuous image of $\mathbf{C}(\aleph_0 \times \kappa) = \mathbf{C}(\kappa)$. Indeed, since $\mathbf{C}(\kappa)$ is compact and $\mathbf{I}(\kappa)$ is Hausdorff, $\mathbf{I}(\kappa)$ is a proper image of $\mathbf{C}(\kappa)$ which is, of course, a closed subspace of $\mathbf{I}(\kappa)$. (A *proper map* is a continuous closed surjective function with each fiber compact. Products of proper maps are proper: Bourbaki [1961, Chapter 1, Section 10 No.2, Lemma 2; p119]. If in addition the fibers are Hausdorff such maps are called *perfect*.) It follows that for X any space, $X \times \mathbf{I}(\kappa)$ is a proper image of $X \times \mathbf{C}(\kappa)$ which is itself a closed subspace of $X \times \mathbf{I}(\kappa)$. Therefore:

**1.2.1 Observation.** If $\mathcal{P}$ is a topological property preserved in proper images and inherited by closed subspaces, then for $\kappa$ infinite $X \times \mathbf{C}(\kappa)$ has $\mathcal{P}$ iff $X \times \mathbf{I}(\kappa)$ has $\mathcal{P}$. In particular, normality of $X \times \mathbf{C}(\kappa)$ is equivalent to that of $X \times \mathbf{I}(\kappa)$

---

2  This cube appears without name or historical reference in Engelking [1989, 2.3I] but not in [1968] and I have not found an earlier reference. Given Engelking's care with citation, I believe it to be his creation.
3  I have found no name in the literature for this cube and, in view of the importance of Arthur's discovery that it is not normal, I have elected to call it the Stone cube.
4  In fact, for any absolute retract $\mathbf{M}$ of a metrizable space of weight $\kappa$ the cube $\mathbf{M}^\omega$ is universal for the class of metrizable spaces of weight not exceeding $\kappa$.

We know, by universality, that spaces can be embedded in these cubes, but it will be helpful to consider some specific embeddings. I do so for $\downarrow\kappa$, $\downarrow\bar{\kappa}$, $\uparrow\kappa$, $\uparrow\bar{\kappa}$, $\kappa$, and $\bar{\kappa}$.

**1.2.2 Observation.** Let $s_\alpha = (s_{\alpha\beta})$ and $t_\alpha = (t_{\alpha\beta})$ where $s_{\alpha\beta} = (\alpha > \beta)$ and $t_{\alpha\beta} = (\alpha \leq \beta)$, and set $S = \{s_\alpha: \alpha \in \kappa\}$, $S^* = \{s_\alpha: \alpha \in \bar{\kappa}\}$, $T = \{t_\alpha: \alpha \in \kappa\}$, and $T^* = \{t_\alpha: \alpha \in \bar{\kappa}\}$. Then:
a) as subspaces of $\mathbf{A}(\kappa)$, S and S* are homeomorphic to $\downarrow\kappa$, and $\downarrow\bar{\kappa}$ respectively;
b) as subspaces of $\mathbf{A}(\kappa)$, T and T* are homeomorphic to $\uparrow\kappa$, and $\uparrow\bar{\kappa}$, respectively, and
c) as subspaces of $\mathbf{C}(\kappa)$, S and S* are homeomorphic to $\kappa$ and $\bar{\kappa}$ respectively.

**Proof.** For $U_\beta = \pi_\beta^{-1}(\{0\})$ open in $\mathbf{A}(\kappa)$ and $\mathbf{C}(\kappa)$, and $V_\beta = \pi_\beta^{-1}(\{1\})$ open in $\mathbf{C}$,
- $U_\beta \cap S = \{s_\gamma : s_{\gamma\beta} = 0\} = \{s_\gamma : \gamma \leq \beta\}$, which corresponds to the interval $[0, \beta+1)$;
- $V_\beta \cap S = \{s_\gamma : s_{\gamma\beta} = 1\} = \{s_\gamma : \gamma > \beta\}$, which corresponds to $(\beta, \kappa)$;
- $U_\beta \cap T = V_\beta \cap S$ while $V_\beta \cap T = U_\beta \cap S$; and
- $V_\beta \cap S^*$ and $U_\beta \cap T^*$ correspond to $(\beta, \kappa]$.

Since S* is compact it is closed in $\mathbf{C}(\kappa)$, a fact easy to observe directly: if x is not in S* there exist $\alpha < \beta$ with $x_\alpha = 0$ and $x_\beta = 1$; $U_\alpha \times V_\beta$ is then a neighborhood of x and its intersection with S* corresponds to $[0,\alpha] \cap [\beta,k] = \varnothing$. Consequently, normality of $X \times \mathbf{C}(\kappa)$ (or of $X \times \mathbf{I}(\kappa)$ ) implies normality of $X \times \bar{\kappa}$. We continue the discussion of X x Y for Y related to $\bar{\kappa}$ in the next subsection and will return to it in 1.5.5, Section 4, and Section 6.

## 1.3 The staircase constructed using $\kappa^+$

A standard construction used to demonstrate that the product of a space X with some space Y is not normal adds a point p to X forming a set X* and creates from that set a space in such a way that the closure of the diagonal of X x X* cannot be separated from the "lid" $X \times \{\bar{p}\}$. For more on this see 6.3. For more general Y there is a frequently used alternative construction which produces a closed set F which one can think of as "below the diagonal" which cannot be separated from some slice $X \times \{\bar{y}\}$ of X x Y. Here I give the case where, as a set, Y equals $\kappa^+$. This construction will be used in 1.5.2, 3.3, and 4.1(b) and referenced in Section 6.

Given a space X which is not compact let $\mathcal{O}$ be an open cover which has no finite subcover, let $\kappa$ be the smallest cardinality of a subcover of $\mathcal{O}$, and let $\{O_\alpha: \alpha \in \kappa\}$ be such a subcover. For each $\alpha$ set $U_\alpha = \bigcup_{\beta < \alpha} O_\beta$ ( with $U_0 = U_1$ ) and notice that for limit ordinals $\alpha$, $\bigcup_{\beta < \alpha} U_\beta = U_\alpha$. (Alternatively, we can just start with an increasing open cover $\{U_\alpha\}$ which satisfies that condition.)

Set $F_\alpha = X \setminus U_\alpha$ and in X x Y set $F = \bigcup_\alpha (F_\alpha \times [0,\alpha])$ and $H = X \times \{\kappa\}$. One can visualize F as a staircase ascending toward H. F is closed in $\uparrow\bar{\kappa}$: suppose $(x,\gamma)$ is not in F; since $\{U_\alpha\}$ is a cover $\bigcap\{F_\alpha\}$ is empty so there exists an index, and therefore a least index $\beta$ such that $x \notin F_\beta$. Note that $\beta$ cannot be a limit ordinal, since for such ordinals $F_\beta = \bigcap_{\alpha < \beta} F_\alpha$ and x would fail to be in some such $F_\alpha$. Hence $(\beta-1,\kappa]$ is open in $\uparrow\overline{\kappa}$. Furthermore, $x \in F_\alpha$ for all $\alpha$ in $[0,\beta)$, $\gamma$ must be in $[\beta,\kappa]$. Thus $U_\beta \times [\beta,\kappa]$ is an open neighborhood of $(x,\gamma)$ which does not meet F. (In closed form, $F = X \setminus \bigcup_\alpha ( U_\alpha \times (\alpha,\kappa])$.)

### 1.4 Some examples of vacuously normal spaces

Spaces with upper or lower topologies, including $\uparrow\kappa$ and $\uparrow\overline{\kappa}$, are vacuously normal. Dowker [1951] used the lower topology $\{\emptyset, \mathbf{R}\} \cup \{(-\infty, r)\}$ on the reals (which in this case coincides with the Scott topology) as an example of a normal space which is not countably paracompact. Steen & Seebach [1970] give several additional examples / constructions of vacuously normal spaces:

| Example(s) | The Set | Topology |
|---|---|---|
| 4 | any set | indiscreet |
| 13-16 | $X^* = X \cup \{p\}$ | $\{X^*\} \cup \{U \subseteq X: U \text{ is open in } X\}$ |
| 52 | $X = (0,1)$ | $\{X\} \cup \{(1/n,1): n>0\}$ |
| 57 | $X = \{2,3,...\}$ | $[ \{\emptyset, X\} \cup \{U_n = \{m: m \text{ divides } n\}: n \geq 2\} ]$ |

In Example 57, [ ... ] indicates the topology generated by the sets listed within.

For X any space the space $X^*$ of examples 13 through 16 is compact (the only open set containing p is X) and vacuously normal (each closed nonempty subset contains p). Thus any space can be embedded as a dense open subspace of a compact vacuously normal space. Note that the topological sum of finitely many nonempty vacuously normal spaces, while not vacuously normal, is normal, and the class of such sums is closed under finite products.

### 1.5 Retracts and very open maps

As is well known, each retraction (continuous map onto a subspace for which the restriction to its range is the identity) preserves normality. In cases involving vacuously normal spaces some results in the converse direction obtain:

**1.5.1 Proposition.** Suppose $r: X \to Z \subseteq X$ is a retraction.
(a) If X is normal, Z is normal.
(b) If Z is normal, the fibers of r are vacuously normal, and r is a closed map, then X is normal.
(c) If X and Z are normal with Z $T_1$ and the fibers of r are vacuously normal, then r is closed.

**Proof.** (a) For F and H disjoint closed subsets of Z set $F^* = r^{-1}(F)$ and $H^* = r^{-1}(H)$. Then $F = F^* \cap Z$ and $H = H^* \cap Z$, so open sets which separate $F^*$ and $H^*$ separate F and H.

(b) Suppose F and H are disjoint closed subsets of X. The intersections of F and H with any fiber are disjoint closed sets, one of which must, by vacuous normality, be empty. Thus r(F) and r(H) are disjoint closed sets, hence separated by disjoint open sets U and V. The disjoint open sets $r^{-1}(U)$ and $r^{-1}(V)$ separate F and H.

(c) If F is a closed subset of X and z is not in r(F), then $H = r^{-1}(z)$ and F are disjoint closed subsets of X and there exists a closed neighborhood N of H which does not meet F. Set $N' = N \cap Z$ and note that for y in N', $r^{-1}(y)$ cannot meet F since it is vacuously normal and meets N. Thus N' is a neighborhood of z which does not meet r(F) and therefore r(F) is closed.

**1.5.2 Example.** There exists a space X with a retraction r: X→Z such that X and Z are normal and each fiber of r is vacuously normal but r not closed.

**Proof.** For $Y = \uparrow\omega$, $Z = \downarrow\overline{\omega}$, and $X = Y \times Z$, let r be the projection from X to Z; Y, Z, and the fibers of r are vacuously normal. The product X is normal (indeed, as we will see in Section 3, vacuously normal) but the staircase $F = \bigcup \{ [n,\omega) \times [0,n]: n \in \omega \}$ is a closed subset of X for which r(F) is not closed.

Factors of a normal $T_1$ product are easily seen to be normal since they are embedded as closed subspaces – for example, $X \times \{y\}$ in $X \times Y$. That proof is not available if $\{y\}$ is not closed, but $X \times \{y\}$ is always a retract so an immediate consequence of 1.5.1(a) is that each factor of a normal product is normal. We will apply 1.5.1(b) to products in Section 3; here I consider another class of maps, suggested by retracts, which also preserve normality and provide a convenient mechanism for generating normal images of product spaces. A retraction from X to Z preserves normality because disjoint open sets of X provide corresponding disjoint open sets of Z. Call a map *very open* if it in fact maps disjoint open sets to disjoint open sets (with the empty set included, so each very open map is open). Homeomorhisms are of course very open and the concepts coincide for maps with Hausdorff domain.

Recall that a space is *irreducible* if its nonempty open sets form a filter base. The definition suggests a wealth of examples - start with your favorite filter base; co-finite subsets of an infinite set and co-countable subsets of an uncountable set are popular. An irreducible space is the closure of each point in the intersection of its filter of open sets - if there are any. Any normal irreducible space (in particular, the closure of any point in a normal space) is vacuously normal (since it has no disjoint nonempty open subsets).

**1.5.3 Proposition.**
(a) If each fiber of an open map is irreducible, that map is very open.
(b) A product of very open maps is very open.
(c) A continuous very open map preserves normality.

**Proof.** (a). At most one of two disjoint open sets can meet an irreducible fiber, so their images are also disjoint.

(b). It suffices to consider basic open sets and thus the product of two maps. For p and q very open with $U_0$ x $V_0$ and $U_1$ x $V_1$ disjoint and open, either $U_0 \cap U_1 = \emptyset$ or $V_0 \cap V_1 = \emptyset$, say the former, *i.e.*, $p(U_0) \cap p(U_1) = \emptyset$. Then $p(U_0)$ x $q(V_0)$ cannot meet $p(U_1)$ x $(V_1)$, as desired.

(c). For q continuous and very open, X normal, and F and H disjoint closed subsets of q(X), $q^{-1}(F)$ and $q^{-1}(H)$ are closed and disjoint, hence contained in disjoint open sets U and V for which q(U) and q(V) are disjoint open sets which separate F and H.

**1.5.4 Example.** Let $\tau_0$ be the functor which carries an arbitrary space to the $T_0$ quotient formed by identifying those points with identical neighborhood systems. Clearly such quotients are very open, so $\tau_0$ carries normal spaces to normal spaces and preserves normality of products. Thus for most considerations of normal products we can, should we wish to do so, restrict our attention to $T_0$ spaces without consequential loss of generality.

**1.5.5 Proposition.** If $X \times \Pi_{\alpha \in A} Y_\alpha$ is normal where $\kappa$ of the $Y_\alpha$ are <u>not</u> vacuously normal then $X \times \mathbf{C}(\kappa)$ and $X \times \bar{\kappa}$ are normal.

**Proof.** Since $Y_\alpha$ is not vacuously normal we can, for each $\alpha$, choose points $s_\alpha$ and $t_\alpha$ such that $\bar{s}_\alpha \cap \bar{t}_\alpha$ is empty. Set $S_\alpha = \bar{s}_\alpha \cup \bar{t}_\alpha$, and let $q_\alpha : S_\alpha \to Z_\alpha = \{0,1\} = \mathbf{2}$ carry $\bar{s}_\alpha$ to $\{0\}$ and $\bar{t}_\alpha$ to $\{1\}$. Each $q_\alpha$ is continuous and very open so for S the product of the $S_\alpha$ and q the product of the $q_\alpha$, q is very open and maps S to $\mathbf{C}(\kappa)$. It follows that $X \times \mathbf{C}(\kappa)$, the image of $X \times S$ under the continuous very open surjective map $\mathrm{id}_X \times q$, is normal, as is its closed subspace $X \times \bar{\kappa}$ is normal.

The corresponding argument with $\kappa$ factors vacuously normal but not indiscrete leads to the conclusion that $X \times \mathbf{A}(\kappa)$ is normal, but we will see stronger results in 3.2. I will return to the subject of normality of products with $\kappa^+$ in Section 4 where it will be shown that 1.5.5 leads, by two different paths, to results such as "If X is not vacuously normal and has all powers normal, then X is compact". Each path uses pseudocompactness. One combines it with closed projections, the other with κ-paracompactness.

## 2 Ancillary properties: pseudocompactness and κ-paracompactness

There is an intimate relationship between normality of uncountable products and pseudocompactness (ranging up to countable compactness) identified by Stone's theorem: **N(Ω)** is not normal. Arthur Stone [1948] established this fact via an elegant, elementary proof by exhibiting two disjoint closed subsets of **N(Ω)** = **N**$^Ω$ which cannot be separated by disjoint open subsets. From that he concluded that if a $T_1$ product is normal, all but countably many of its factors must be countably compact (because each $T_1$ space which is not countably compact contains a copy of **N** as a closed subspace). That conclusion may not hold for spaces which are not $T_1$, but we can conclude that countably many factors fail to have a property which we will call "loose compactness" which for our purposes is strong enough.

### 2.1 Loosely compact spaces

Our interest is in conditions stronger than *pseudocompactness* (each continuous real-valued function bounded) but weaker than countable compactness and consider two such conditions here. Recall (*e.g.*, Gillman & Jerison [1960, 1.10-1.15] ) that a subset is a < zero set, cozero set > if it equals < $f^{-1}(\{0\})$, $f^{-1}(\{$**R**$-\{0\}\})$ > for some continuous real-valued function f and that a collection of subsets is < locally finite, discrete > if each point has a neighborhood which meets < finitely many, at most one > of the members. A space is pseudocompact if and only if each locally finite family of nonempty cozero sets is finite (for $\{C_n\}$ such a family and, for each n, $f_n$ a continuous real-valued function which is zero off of $C_n$ and at least 1 somewhere on $C_n$, the function $\Sigma_n n f_n$ is continuous since $\{C_n\}$ is locally finite and bounded if and only if $\{C_n\}$ is finite).

A space is *lightly compact* - Bagley, Connell, & MacKnight [1958] - aka "feebly compact" - Scarborough & Stone [1966] - if each locally finite family of nonempty open subsets is finite or, equivalently (Stone [1960]), if each locally finite family of disjoint open sets is finite. Other characterizations can be found in McCoy [1973] and Özer [1983]. Call a space *loosely compact* if each discrete family of nonempty closed subsets is finite. Notice that discrete sets must be disjoint so each vacuously normal space is, almost trivially, loosely compact.

**2.1.1 Observation.** Of the four properties of a space X which follow, (a) ⇨ (b) ⇨ (d) and (a) ⇨ (c) ⇨ (d) but no other implications hold in general. If X is completely regular, (d) ⇨ (c); if X is normal (d) ⇨ (b); and if X is normal $T_1$, (d) ⇨ (a) and (a) through (d) are equivalent. In particular, (c) ⇨ (b) for X normal, while in the other direction (b) ⇨ (c) if X is regular or $T_1$.

(a) countable compactness;
(b) loose compactness;
(c) light compactness;
(d) pseudocompactness.

**Proof.** (a) ⇨ (b): Suppose $\{S_n\}$ is a discrete family of closed subsets and, for each n, set $F_n = \bigcup_{m \geq n} S_m$. Then $\bigcap_n F_n$ is empty and, since $\{S_n\}$ is discrete, each $F_n$ is closed; hence by countable compactness some $F_n$ is empty. Thus all but finitely many of the $S_n$ are empty.

(b) or (c) ⇨ (d): If X were not pseudocompact there would exist an unbounded continuous function f: $X \to \mathbf{R}$ with, we may assume, $\mathbf{N}$ contained in the range of f. But then $\{f^{-1}((n))\}$ would be an infinite discrete family of nonempty closed subsets, contrary to (b) and $\{f^{-1}((n-1/3, n+1/3))\}$ would be an infinite locally finite family of nonempty open subsets, contrary to (c). That no other implications hold in general will be shown by the examples of 2.1.2.

(d) ⇨ (c) for X completely regular: suppose $\{U_n\}$ is a locally finite family of nonempty open subsets of X and, for each n let $f_n: X \to [0,1]$ be continuous with $f_n(x_n) = 1$ for some point $x_n \in U_n$ and $f_n(X \setminus U_n) = 0$ (complete regularity ensures the existence of such functions). Set $f = \Sigma_n n f_n$; since $\{U_n\}$ is locally finite f is continuous and since X is pseudocompact it must be bounded, hence $\{U_n\}$ is finite.

(d) ⇨ (b) for X normal: suppose $\{S_n\}$ is a discrete family of nonempty closed sets of X. Then $S = \bigcup_n S_n$ and each $T_n = \bigcup_{m \neq n} S_m$ is closed, hence each $S_n = S \setminus T_n$ is both open and closed in S. As a closed subspace of a normal pseudocompact space S is pseudocompact (*e.g.*, by the Tietez extension theorem), so $\{S_n\}$ must be finite, as desired.

(d) ⇨ (a), X normal $T_1$: if the Hausdorff space X is not countably compact it contains an infinite closed discrete subspace $\{s_n\}$ and, by normality, the continuous function $f(s_n) = n$ extends to X, so X is not pseudocompact.

(b) ⇨ (c) for X regular or $T_1$: using Arthur Stone's characterization of light compactness, suppose $\{U_n\}$ is an infinite locally finite family of <u>disjoint</u> nonempty open subsets and choose $x_n \in U_n$. Since $\{U_n\}$ are disjoint, $\{\overline{x}_n\}$ is an infinite family of closed sets. If X is regular, each $x_n$ has a closed neighborhood contained in $U_n$ so $\overline{x}_n$ is contained in $U_n$, as is of course true if X is $T_1$. Thus $\{\overline{x}_n\}$ is a locally finite family of disjoint closed sets. If each neighborhood of y meets $\bigcup \{\overline{x}_n\}$, then by local finiteness it must be the case that each neighborhood of y meets some particular $\overline{x}_n$, and hence that y is in $\overline{x}_n$. Thus either y has a neighborhood which meets no $\overline{x}_n$, or y is in some $\overline{x}_n$, and $U_n$ is a neighborhood of y which meets no other. Thus $\{\overline{x}_n\}$ is discrete.

Except as it pertains to loosely compact spaces, Observation 2.1.1 is included in Theorems 3 and 4 of to Bagley, Connell, & MacKnight [1958]; the portions relating countable compactness and pseudocompactness are due to Hewitt [1948, Theorem 30] who, at [1948, page 69], noted (as discussed below) that the deleted Tychonoff Plank is lightly compact but not countably compact. I should note that Observation 2.1.1 is far from comprehensive: Stephenson Jr., [2004] considers countable compactness, light compactness, and pseudocompactness plus five additional conditions which lie between light compactness and pseudocompactness; Steen & Seebach define a "weak countable compactness" which does not imply pseudocompactness; and variations of pseudocompactness which strengthen or weaken it abound: see, *e.g.*, Garcia-Ferreira & Ortiz-Castillo [2018] and Dorantes-Aldama, Okunev, & Tamariz-Mascarúa [2018].

**2.1.2 Examples.** There exist spaces which are:
(a) not countably compact but lightly compact and loosely compact;
(b) loosely compact but not lightly compact;
(c) lightly compact but not loosely compact;
(d) neither loosely compact nor lightly compact but pseudocompact.

**Proof.** (a). The space ↓**N** (N with the lower topology $\{N\} \cup \{[0,n): n \in N\}$ ) is vacuously normal, hence loosely compact; lightly compact; but not countably compact.

(b). Let $X \subset N \times N$ be $\{(m,n): n \leq m\}$ with $\{\varnothing\} \cup \{\Delta_{mn}: n \leq m\}$ forming a base for the open sets where $\Delta_{mn}$ is the triangle with apex at (m,n) spreading to a base from (m-n,0) to (m+n,0). The one point sets $\Delta_{m0} = \{(m,0)\}$ form an infinite locally finite family of open subsets, so X is not lightly compact, but since each point (m,n) has, in its closure, all points (i,i) for $i \geq m$, each pair of nonempty closed sets meet so X is vacuously normal, hence loosely compact.

(c). The deleted Tychonoff plank $\Gamma = \Omega^* \times N^* \setminus \{(\Omega,\omega)\}$ is pseudocompact Tychonoff, hence lightly compact, but not loosely compact: $\{(\Omega,n)\}$ is closed and discrete.

(d). Bagley, Connell, & MacKnight [1958, p502] give an example of a space which is pseudocompact $T_1$ but not lightly compact (hence, being $T_1$, not loosely compact). As an example motivated by theirs, let X be two copies of **N** x **N**, X = **N** x **N** x **II,** where the first copy, $X_0 = $ **N** x **N** x {0} is open and carries its usual discrete topology but points of the second copy, $X_1 = $ **N** x **N** x {1} have a coarser topology in which each neighborhood of (m,n,1) contains not only (m.n,0), but also those points in $X_0$ above or to the right of it. That is, for $_{mn}L = \{(i,j): i = m$ and $j \geq n$ or $i \geq m$ and $j = n\}$, each neighborhood of (m,n,1) includes itself plus $_{mn}L$ x {0}.

Notice that the closure of (m,n,0) is {(m,n,0), (m,n,1)} plus those points of $X_1$ to the left of or below (m,n,1). Thus the closure of $_{mn}L$ includes both the row and column of $X_1$ through (m,n,1) (and the corresponding row and column in $X_0$ plus much more). Since $_{mn}L$ is the smallest neighborhood of (m,n,1) and each real-valued continuous function is constant on the closure of the smallest neighborhood of a point (or the closure of the point if no such neighborhood exists), each such function on X must be constant on a row plus column of $X_1$, indeed, each row and column of X. Thus X is pseudocompact. The diagonal {(n,n,0)} of $X_0$ is an infinite locally finite family of open sets so X is not lightly compact, and the points of $X_1$ are closed and discrete so X is not loosely compact.

**2.1.3 Theorem.** If X x Y is loosely compact, its projections are z-closed (map zero sets to closed sets).

**Proof.** Suppose $\pi = \pi_Y$ is not z-closed and let f: X x Y → **I** be a continuous function such that for $Z = f^{-1}(0)$, $\exists z \in \overline{\pi Z} \setminus \pi Z$. Replacing f with min{ max{ 0, f(x,y)/f(x,z) }, 1 }, we may assume that $f \in C(X \times Y, \mathbf{I})$ and f(X x {z}) = {1}. For each x in X set $V(x) = \pi((\{x\} \times Y) \setminus Z)$; each V(x) is a neighborhood of z such that {x} x V(x) does not meet Z.

We will constructs points $(x_n, y_n) \in Z$ with $y_{n+1} \in \bigcap_{m \leq n} V(x_m)$ and $f(x_n, y_{n+1}) > n/n+1$. To start, choose $(x_0, y_0) \in Z$ and suppose, inductively, that for m ≤ n satisfactory points $(x_m, y_m)$ have been found. Set $V = (\bigcap_{m \leq n} V(x_m))$ and $W = \pi(f^{-1}((n/n+1, 1]) \cap (\{x_n\} \times Y))$. Since both V and W are neighborhoods of z and $z \in \overline{\pi Z}$, X x (V ∩ W) must meet Z, so we can choose $(x_{n+1}, y_{n+1})$ from that intersection. Then $(x_{n+1}, y_{n+1}) \in Z$, $y_{n+1} \in \bigcap_{m \leq n} V(x_m) = V$, and since $y_{n+1} \in W$, $f(x_n, y_{n+1}) > n/n+1$, as desired.

Let $F_n$, $H_n$, and $S_n$ be the closures of $x_n$, $y_n$, and $(x_n, y_n)$ respectively, and note that $S_n = F_n \times H_n$. Thus if $x \in F_n$, f(x,y) equals $f(x_n, y)$ for any y, and $x \in F_m \cap F_n$ where m < n would imply both $f(x, y_n) = f(x_m, y_n) > 0$ since $y_n \in V(x_m)$ and $f(x, y_n) = f(x_n, y_n) = 0$ because $(x_n, y_n)$ is in Z. Therefore the $F_n$ are disjoint and consequently the $S_n$ are disjoint.

As X x Y is loosely compact, {$S_n$} has a cluster point, (p,q) where, since $S_n = F_n \times H_n$, p is a cluster point of {$F_n$} and q is as cluster point of {$H_n$}. But this implies both $f(p,q) = \lim_n f((x_n, y_n)) = 0$ and $f(p,q) = \lim_n f((x_n, y_{n+1})) = 1$, which contradiction establishes that $\pi_Y$ is z-closed.

The proof of 2.1.3 is adapted from a proof, suggested in Hager [1969, Theorem 4.2] and provided as Comfort & Hager [1971, Theorem 4.1]. As Comfort and Hager acknowledge, the proof itself owes a debt to proofs of various related properties by Frolík, Glicksberg, and Isbell.

**2.1.4 Corollary.** If X x Y is pseudocompact normal and Y is $T_1$, then $\pi_Y$ is closed.

**Proof.** By 2.1.1 (d) $\Rightarrow$ (b) X x Y is loosely compact hence by 2.1.3 $\pi_Y$ is z-closed. It follows that $\pi_Y$ is closed since any continuous z-closed map p with normal domain and $T_1$ range is closed: given disjoint closed sets F and $p^{-1}(y)$ there exists, by Urysohn's Lemma, a zero set Z containing F disjoint from $p^{-1}(y)$ and p(Z) is a closed set containing p(F) which does not contain y.

**2.1.5 Proposition.** If $X \times \Pi_{\alpha \in A} Y_\alpha$ is normal and no $Y_\alpha$ is loosely compact, then A is countable.

**Proof.** Since each factor $Y_\alpha$ is not loosely compact it contains an infinite discrete family of nonempty closed subsets $\{S_{\alpha n}\}$ with, we may assume, each $S_{\alpha n}$ irreducible. The quotient $q_\alpha(S_{\alpha n}) = \{n\}$ from $S_\alpha = \bigcup_n S_{\alpha n}$ to **N** is very open and since $\{S_{\alpha n}\}$ is discrete, $S_\alpha$ is closed. Hence for $S = \Pi_\alpha S_\alpha$ and $q = \Pi_\alpha q_\alpha$, q is very open and continuous with range $\mathbf{N}^\kappa$, where $\kappa = |A|$. Since $\Pi_{\alpha \in A} Y_\alpha$ is normal, S and thus $\mathbf{N}^\kappa$ must be normal, so by Stone's theorem the cardinal $\kappa$ must be countable.

Referring again to the case $X^\kappa$ normal for $\kappa$ uncountable, 2.1.5 implies that if X is not vacuously normal then $X^\kappa$ is loosely compact, from which it follows that $X^\kappa \times \overline{\kappa}$ is pseudocompact (as is any product of a pseudocompact space with a compact space). By 1.5.5 the product $X^\kappa \times \overline{\kappa}$ is also normal so by 2.1.4 its projections are closed which implies that $X^\kappa$ is $[\omega,\kappa]$-compact. (As will be discussed in Section 4.)

## 2.2 κ-paracompact spaces

A space is ***κ-paracompact*** if each κ-fold open cover has a locally finite open refinement. As with the relationship between normality of uncountable products and pseudocompactness, there is an intimate relationship between normality of finite or countable products and κ-paracompactness. For finite products that relationship is suggested by Dowker's [1951, Theorem 4] proof that for X normal, X x **I** is normal if and only if X is countably paracompact and the extension of Dowker's result by Morita [1962, Theorem 2.4]: X x **I**(κ) is normal if and only if X is κ-paracompact. For countable products it is exposed in Zenor [1971, Theorem A], announced in [1969]: A countable product of non-trivial normal $T_1$ spaces is normal if and only if it is countably compact.

The results in this subsection are mostly known. I review them here because they are often presented without proof, in the presence of extraneous separation assumptions, and/or as consequences of less transparent proofs of more general results.

By listing a collection as, for example, $\{O_\alpha: \alpha \in \kappa\}$, I indicate only that it is ***κ-fold*** - has κ or fewer members (empty sets and repetitions allowed). Such a collection is ***directed*** if each pair of members (and thus any finite sub-collection of members) is contained in a member and ***increasing*** ( aka "monotone increasing" or "monotone" ) if in fact $\alpha < \beta$ implies $O_\alpha \subseteq O_\beta$. An increasing family $\{O_\alpha: \alpha \in \kappa\}$ is ***well indexed*** if for each β $\bigcup_{-1 \leq \alpha < \beta} O_{\alpha+1}$ equals $O_\beta$. Given a collection $\{O_\alpha\}$, the associated well indexed collection $\{U_\alpha = \bigcup_{-1 \leq \alpha < \beta} O_{\alpha+1}\}$ has the same union and the same or smaller cardinality as $\{O_\alpha\}$. We have already seen a use of well indexing in the staircase examples.

Where $\mathcal{O} = \{O_\alpha: \alpha \in A\}$ and $\mathcal{U} = \{U_\alpha: \alpha \in B\}$ with $B \subseteq A$ and each member of $\mathcal{U}$ is contained in a member of $\mathcal{O}$ ( in particular where $\mathcal{U}$ is a ***refinement*** - has in addition the same union ), the collection $\mathcal{U}$ is ***index faithful*** or ***faithfully indexed*** ( in relation to $\mathcal{O}$) if for each $\alpha \in B$ $V_\alpha \subseteq U_\alpha$. A ***shrinking*** of an open cover $\{U_\alpha\}$ is a faithfully indexed open refinement $\{V_\alpha\}$ where, in fact, $\overline{V}_\alpha \subseteq U_\alpha$. for each α. The DeMorgan dual of a shrinking, a family of closed sets with empty intersection enclosed by open sets whose closures also have empty intersection, is called an ***expansion***.

### 2.2.1 Remarks.
(a) An open cover has an open locally finite refinement if and only if it has a faithfully well indexed open locally finite refinement.
(b) An open cover has an open locally finite refinement by sets whose closures also form a refinement if and only if it has a locally finite shrinking.

**Proof.** The "if" directions are trivial. Let $\mathcal{O} = \{O_\alpha: \alpha \in \kappa\}$ be a κ-fold open cover with $\mathcal{U}$ an open locally finite refinement of $\mathcal{O}$.

(a): for each α set $V_\alpha = \bigcup \{U \in \mathcal{U} : U \subseteq O_\alpha \text{ but } U \not\subseteq O_\beta \text{ for } \beta < \alpha\}$. Since $\mathcal{U}$ is a refinement each of its members is contained in some first $O_\alpha$ so $\{V_\alpha\}$ is a cover, and since $V_\alpha \subseteq O_\alpha$ it is a faithfully indexed refinement of $\mathcal{O}$. Furthermore, $\{V_\alpha\}$ is locally finite: a neighborhood which meets $V_\alpha$ meets a member $U_\alpha$ of $\mathcal{U}$ used in the construction of $V_\alpha$ and $U_\alpha$ is necessarily distinct from such a $U_\beta$ for $\beta < \alpha$. For $W_\alpha = \{\bigcup_{-1 \leq \beta < \alpha} V_{\beta+1}\}$, $\{W_\alpha\}$ is a faithfully well indexed open locally finite refinement of $\mathcal{O}$.

(b): Suppose that, for $\mathcal{U}$ as above, each U in $\mathcal{U}$ has $\overline{U}$ contained in some member of $\mathcal{O}$. For each α set $V_\alpha = \bigcup \{U \in \mathcal{U} : \overline{U} \subseteq O_\alpha \text{ but } \overline{U} \not\subseteq O_\beta \text{ for } \beta < \alpha\}$. Again, each $\overline{U}$ is contained in some first $O_\alpha$ so $\{V_\alpha\}$ is a faithfully indexed refinement of $\mathcal{O}$. Since $\mathcal{U}$ is locally finite $\overline{\mathcal{U}}$ is locally finite, so $\overline{V}_\alpha \subseteq \text{cl}(\bigcup \{U \in \mathcal{U} : U \subseteq O_\alpha\}) = \bigcup \{\overline{U} \in \mathcal{U} : \overline{U} \subseteq O_\alpha\} \subseteq O_\alpha$, hence $\{V_\alpha\}$ is a shrinking, as desired.

**2.2.2 Proposition.** For a space X the conditions which follow are equivalent:
(a) X is $\kappa$-paracompact;
(b) each directed $\kappa$-fold open cover of X has an open locally finite refinement;
(c) each directed $\kappa$-fold open cover of X has a locally finite shrinking;
(d) each well indexed $\kappa$-fold open cover of X has a locally finite shrinking.

**Proof.** (a) $\Rightarrow$ (b) as a special case. (b) $\Rightarrow$ (c): let $\mathcal{O} = \{O_\alpha: \alpha \in A\}$ be an $\kappa$-fold directed open cover of X and let $\mathcal{U} = \{U_\alpha: \alpha \in A\}$ be a faithfully indexed open locally finite refinement. For $B \subseteq A$ finite non-empty and $\beta \in B$ set $F_{B(\beta)} = \overline{U}_\beta \setminus \bigcup_{\alpha \in A \setminus B} U_\alpha$. Since $\{U_\alpha\}$ is locally finite and a cover, each point has a neighborhood which meets only finitely many of the $U_\alpha$ and is contained in one of them, that is, a neighborhood contained in some $F_{B(\beta)}$. Hence the interiors of $\{F_{B(\beta)}\}$ form a cover. Also, since $F_{B(\beta)} \subseteq \bigcup_{\alpha \in B} U_\alpha \subseteq \bigcup_{\alpha \in B} O_\alpha$ and $\mathcal{O}$ is directed, $\{F_{B(\beta)}\}$ refine $\mathcal{O}$. Let $\mathcal{V}$ be the collection of finite unions of interiors of members of $\{F_{B(\beta)}\}$; then both $\mathcal{V}$ and $\overline{\mathcal{V}}$ refine $\mathcal{O}$ and $\mathcal{V}$ is directed, so it has an open locally finite refinement $\mathcal{W}$, and $\overline{\mathcal{W}}$ also refines $\mathcal{O}$. Thus ( by 2.2.1(a) ) $\mathcal{O}$ has a locally finite shrinking.

(c) $\Rightarrow$ (d) is another special cases. (d) $\Rightarrow$ (a). Proceeding by induction, suppose $\kappa$ is the smallest cardinal for which the assertion fails, suppose X has the indicated property, let $\{O_\alpha: \alpha \in \kappa\}$ be a $\kappa$-fold open cover of X, let $\{U_\alpha: \alpha \in \kappa\}$ be the associated well indexed open cover $U_0 = \varnothing$, $U_\alpha = \bigcup_{\beta < \alpha} O_\beta$, and let $\{V_\alpha\}$ be a $\kappa$-fold locally finite shrinking of $\{U_\alpha\}$, so $\overline{V}_\alpha \subseteq U_\alpha$. The closed subspace $\overline{V}_\alpha$ of X satisfies the hypotheses on X, hence by induction there exists an $\alpha$-fold locally finite family $\mathcal{U}_\alpha$ of sets open in $\overline{V}_\alpha$ which refine the trace on $\overline{V}_\alpha$ of $\{O_\beta: \beta \leq \alpha\}$. Let $\mathcal{V}_\alpha = \{U \cap V_\alpha: U \in \mathcal{U}_\alpha\}$ and $\mathcal{V} = \bigcup_\alpha \mathcal{V}_\alpha$. Since $\{V_\alpha\}$ and each $\mathcal{V}_\alpha$ is locally finite, $\mathcal{V}$ is a locally finite refinement of $\{O_\alpha\}$.

In 2.2.2 some restriction on the covers which shrink is necessary since (as we will see) shrinking of all finite open covers is equivalent to normality and, for example, a countably paracompact space (indeed, a countably compact locally compact Tychonoff space - Steen & Seebach [1970, Example 106] ) need not be normal. Directed covers provide a natural restriction in this case since, except those covers which include the covered space itself, a non-empty directed cover must be infinite. Normality removes that particular reason for a restriction, and indeed, removes every need for it.

**2.2.3 Proposition.**
(a) A space is normal if and only if each 2-fold (or finite) open cover shrinks.
(b) A space is $\kappa$-paracompact normal if and only if each $\kappa$-fold open cover of it has a locally finite shrinking.

**Proof.** (a): For F and G disjoint closed subsets of X, U and V are open sets with $F \subseteq U$ and $H \subseteq V$ iff $\{X \setminus U, X \setminus V\}$ is a closed refinement of the open cover $\{X \setminus F, X \setminus H\}$, and $\overline{U}$ and $\overline{V}$ are disjoint iff $\{X \setminus V, X \setminus U\}$ is a shrinking.

(b): Clearly a space in which each $\kappa$-fold open cover has a locally finite shrinking is $\kappa$-paracompact and, in view of (a), normal. Suppose X is $\kappa$-paracompact normal, let $\mathcal{O} = \{O_\alpha : \alpha \in A\}$ be a $\kappa$-fold open cover of X, let $\mathcal{U}$ be the collection of unions of finite sub-collections of $\mathcal{O}$, and, invoking 2.2.2, let $\mathcal{W}$ be a shrinking of $\mathcal{U}$. For W in $\mathcal{W}$ let B be a minimal subset of A such that $W \subseteq \bigcup_{\alpha \in B} O_\alpha$, say with $|B| = n$. It suffices to show, for $\beta \in B$, that W can be replaced by open sets P and Q with $W \subseteq P \cup Q$, $\overline{P} \subseteq O_\beta$, and $\overline{Q}$ contained in the union of the remaining n-1 members of $\mathcal{O}$ indexed by B, since by repetition of this process $\mathcal{W}$ can be converted to a shrinking of $\mathcal{O}$.

Set $U = O_\beta$ and $V = \bigcup_{\alpha \in B \setminus \{\beta\}} O_\alpha$, and, invoking normality, let T be an open set with $W \subseteq T$ and $\overline{T} \subseteq U \cup V$. Set $F = \overline{T} \setminus U$ and $H = \overline{T} \setminus V$; since the closed subspace $\overline{T}$ is normal there exist open subset R and S of $\overline{T}$ with $F \subseteq R$, $H \subseteq S$ and $\overline{R} \cap \overline{S} = \varnothing$. For $P = T \setminus S$ and $Q = T \setminus R$, $W \subseteq P \cup Q = T$, $\overline{P} \subseteq U = O_\beta$, and $\overline{Q} \subseteq V$, as desired.

**2.2.4 Corollary:** In a normal $\kappa$-paracompact space, each $\kappa$-fold open cover has a locally finite cozero set refinement. Consequentially a pseudocompact normal space is $\kappa$-paracompact if and only if it is $[\omega,\kappa]$-compact.

**Proof.** For $\{O_\alpha\}$ a locally finite refinement, $\{U_\alpha\}$ a shrinking of it, $f_\alpha$ a Urysohn function zero on the compliment of $O_\alpha$ with $f_\alpha(\overline{U}_\alpha) = 1$, and $C_\alpha = f_\alpha^{-1}((0,1])$, $\{C_\alpha\}$ are as desired and, in a pseudocompact space, must (like any locally finite family of cozero sets) be finite.

The equivalences among locally finite refinements and faithfully indexed or well indexed refinements are often used without comment. The content of 2.2.2 appears, for example, in Mack [1967, Theorem 5]. The fact that each open locally finite open cover of a normal space has a closed refinement dates back to Dieudonné [1944, Theorem 6]. The existence of a cozero set refinement is implicit in Michael's proof (for paracompact spaces) of [1953, Proposition 2] and thus so, almost, is the existence of shrinkings. There are suggestions of 2.2.3 in Mack [1967], Nagami [1970], Morita [1975], and Rudin [1978], and shrinkings are considered, along with variations, in Rudin [1985]. Michael [1957] characterizes paracompactness in terms of closure preserving covers, and Katuta [1977] focus on characterizations in terms of increasing covers, and in particular the cushioned covers introduced by Michael [1959]. Other characterizations involve star refinements, barycentric refinements, and countable collections s which together make up a refinement. The fact that a pseudocompact normal space is $\kappa$-paracompact if and only if it is $[\omega,\kappa]$-compact is due to Morita [1962, Theorem 1.8].

Arguably the most important topological result concerning paracompactness is Arthur Stone's [1948] discovery that each metrizable space is paracompact. At its core, Stone's result does not require Hausdorff separation. As Engelking [1987, Remark 4.4.2] notes, the proof (at least, Mary Ellen Rudin's [1969] proof) of Stone's Theorem actually shows:

**2.2.5 Theorem (A. H. Stone):** Each open cover of a pseudometrizable space has an open refinement which is both locally finite and σ-discrete (a union of countably many discrete families).

**2.2.6 Corollary:** Suppose f is a continuous real valued function on a product X x Y and ☉ is a cover of X made up of sets of the form $\pi_X(f^{-1}(O) \cap (X \times \{y\}))$ where O is open in the range of f and y is in Y. If $\pi_X$: X x Y →X z-closed, then ☉ has an open σ-discrete locally finite refinement.
**Proof.** We may suppose the range of f is contained in **I**; notice that, since $\pi_X$ is z-closed, for any continuous function g: X x Y→**I** the function $\sup_Y g(x) = \sup\{g(x,y): y \in Y\}$ is continuous ($\sup_Y g^{-1}([0,s]) = X \setminus \pi_X g^{-1}((s,1])$ and $\sup_Y g^{-1}([r,1]) = \bigcap_n \pi_X g^{-1}([r-1/n,1])$ are both closed from which it follows that each $\sup_Y g^{-1}([r,s])$ is closed). Thus in particular for $g_w(x) = |f(w,y) - f(x,y)|$, $\sup_Y g_w$ is continuous.

For (w,x) in X x X set $\rho(w,x) = \sup_Y g_w(x)$; ρ is a pseudometric as it inherits that property from the absolute value function, and is continuous since the ball of radius ε centered at x, $\sup_Y g_x^{-1}([0,\varepsilon))$, is open. Clearly each restriction of f to slices X x {y} is continuous with respect to ρ, so our conclusion follows from 2.2.5.

Notice that, in fact, f is continuous on X x Y when X carries the topology induced by the pseudometric ρ: for B the ball of radius ε/2 centered at x, y fixed, r = f(x,y), J the intersection of ( r-ε/2,r+ε/2 ) with **I**, and V = $f^{-1}(J)$, B x V is a neighborhood of (x,y) contained in $f^{-1}((r-\varepsilon,r+\varepsilon) \cap \mathbf{I})$. The use of z-closed projections to ensure that $\sup_Y g$ is continuous is reasonably precise: if $\sup_Y g$ is continuous for each g in C(X x Y, **I**), then $\pi_X$ is z-closed. ( For $w \notin f^{-1}(0) = Z$ and $g(x,y) = \min\{f(x,y)/f(w,y),1\}$, $\sup_Y g(w) = 1$ but $\sup_Y g(\pi_X Z) = 0$; Noble [1969b, Theorem 1], due to Comfort & Hager. )

## 3  The implicated properties: vacuously normal and compact normal spaces

Of course, compact Hausdorff spaces are normal and well understood, so for compact normal spaces our focus here will be spaces which are not Hausdorff, and therefore neither $T_1$ nor regular. I start with a characterization of vacuously normal spaces. (Spaces in which each pair of nonempty closed subsets meet.) Call a point x ∈ X a **vn-point** of X if $\overline{x}$ meets $\overline{y}$ for each point y ∈ X.

**3.1 Proposition**
(a) A space is vacuously normal if and only if each of its points is a vn-point.
(b) A point in a product space is a vn-point if and only if each of its factors is a vn-point.
(c) A product is vacuously normal if and only if each of its factors is vacuously normal.
**Proof.** (a) is trivial and (b) follows immediately from the observation that for $x = (x_\alpha)$, its closure $\bar{x}$ equals $\Pi_{\alpha \in A} \bar{x}_\alpha$. Given (a) and (b), (c) is obvious.

**3.2 Proposition.** Suppose X is vacuously normal and Y normal. If $\pi_Y: X \times Y \to Y$ is closed then $X \times Y$ is normal. If Y is $T_1$, the converse holds.
**Proof.** This result is an immediate consequence of 1.5.1 (b) and (c) since $\pi_Y$ is a retraction (under the identification of Y with some $\{x\} \times Y$).

**3.3 Example.** The product of a vacuously normal space with a compact Hausdorff space need not be normal.
**Proof.** In $\downarrow\omega \times \bar{\omega}$ the staircase $F = \bigcup_n \{[n,\infty) \times [0,n]\}$ cannot be separated from $\omega \times \{\infty\}$ since the smallest open set containing F is $\omega \times \omega$.

Turning our attention to compact normal spaces, I start with a class of such spaces which is significant in the non-Hausdorff context because it includes each closure of a point in a compact normal space. Call a space X ***trivially compact*** if each open cover of X includes, as one of its members, the space X itself, *i.e.,* if and only if it contains a point whose only neighborhood is X. Notice that trivial compactness is preserved by products.

**3.4 Observation.** Each irreducible compact normal space is trivially compact and vacuously normal.
**Proof.** As was noted earlier, any irreducible normal space is vacuously normal because it has no nonempty disjoint open subsets. Consequently its nonempty closed subsets form a filter base and if, in addition, the space is compact that filter base has nonempty intersection. The only open set which meets that intersection is the space itself, which therefore must be included in any open cover.

Note that in terms of the order of specialization $x \leq y \Leftrightarrow y \in \bar{x}$ a point x is a least element of $\bar{x}$ while points in the intersection of its nonempty closed subsets are co-equal maximums (an open neighborhood of one excluding others would have, as its compliment, a nonempty closed subset from which it was excluded). Hence if $\bar{x}$ is a compact normal $T_0$ space, co-equal points coincide and the intersection of the closed subsets of $\bar{x}$ is a single point. For example, $\bar{x}$ might be some $\bar{\kappa}$ with the lower topology.

**3.5 Theorem.**
(a) A $T_0$ space X is compact normal if and only if there exist a retraction r: $X \to Z$ where Z is compact Hausdorff and each fiber of r is trivially compact and vacuously normal.
(b) A product is compact normal if and only if its factors are compact normal.

**Proof.** (a). Suppose X is a compact normal $T_0$ space; for x in X let r(x) be the closed point of X which is the intersection of the nonempty closed subsets of $\bar{x}$ and set Z=r(X). Let S(x) be the saturation, in X, of r(x) ( the intersection of the open sets which contain r(x) ). Notice that S(x) is trivially compact (each open cover includes a subset which contains r(x) and thus S(x) ) and vacuously normal (each nonempty closed subset contains r(x) ). Furthermore, $S(x)=r^{-1}(r(x))$: $y \in S(x)$ iff $r(x) \in \bar{y}$ iff $r(x)=r(y)$ since r(y) is the only closed point in $\bar{y}$. Thus the fibers of r are trivially compact and vacuously normal.

To see that r is continuous suppose $U \subseteq X$ is open and $x \in r^{-1}(U \cap Z)$. Then $r(x) \in U$ and $\{r(x)\}$ and $X \setminus U$ are disjoint closed subsets of the normal space X, so U contains a closed neighborhood N of r(x). For y in N, $r(y) \in N$ since $r(y) \in X \setminus N$ would imply that S(y) and hence y were in $X \setminus N$. Thus $r(N) \subseteq U$ and therefore r is continuous.

As a continuous image of a compact space Z is compact and by 1.5.1(a) it is normal hence (since its points are closed) Haudorff. Thus the retraction r and subspace Z are as desired.

For the converse, X is compact since any finite subcover of Z by open subsets of X covers $r^{-1}(Z) = X$ and as a continuous map from a compact space to an Hausdorff space the map r is closed, so by 1.5.1(b) X is normal.

(b). As a retraction, each factor of a compact normal product is compact normal. Since a product of retractions with compact Hausdorff range is itself a retraction with compact Hausdorff range and a product of trivially compact vacuously normal fibers is trivially compact and, by 3.1(c), vacuously normal, it is a consequence of (a) that a product of compact normal $T_0$ spaces is normal (and, of course, compact). It follows, since (as was observation in 1.5.4) the $\tau_0$ functor is very open, that a product of compact normal spaces is compact normal.

**4 Test products: X x Y where as a set $Y = \kappa^+$.**

I return to the consideration of $X \times \bar{\kappa}$ and related spaces begun in Section 1.2. Recall that $\bar{\kappa}$ is the set $\kappa^+ = [0,\kappa]$ with the (compact) order topology and $\uparrow\bar{\kappa}$ is $\kappa^+$ with the upper topology $\{\emptyset, \kappa^+\} \cup \{(\alpha.\kappa]: \alpha \in \kappa\}$. For $\tau$ a topology let $\delta(\tau)$ be the least cardinal such that each intersection of fewer than $\delta(\tau)$ members of $\tau$ is again in $\tau$. (If $\tau$ is Alexandroff discrete, *i.e.*, closed under the formation of intersections, set $\delta(\tau)=\infty$.)

We will use topologies on $\kappa^+$ to "test" cardinality related properties of a companion factor, so call such a topology $\tau$ "test worthy" if it contains the topology of $\uparrow\overline{\kappa}$ and for each cardinal $\lambda$ in $[\delta(\tau),\kappa]$ the interval $[\lambda,\kappa]$ is not open. Note that topologies between $\uparrow\overline{\kappa}$ and $\overline{\kappa}$ are testworthy, as well as non-discrete topologies $\tau$ with $\delta(\tau) = \kappa$ such as the topology, for $\kappa$ regular, with proper subsets open if and only if they do not contain $\kappa$ or have compliment of cardinality less than $\kappa$ - the "co- $< \kappa$" topology. A range of topologies between $\uparrow\overline{\kappa}$ and $\overline{\kappa}$ (or $\downarrow\overline{\kappa}$ and $\overline{\kappa}$) can be generated by considering $\kappa^+$ as a subspace of $\mathbf{A}(\lambda) \times \mathbf{C}(\mu)$ where $\lambda+\mu=\kappa$.

A collection $\mathcal{B}$ of open subsets of a space Y is a *hyperbasis* for Y if for each closed subset $F \subseteq Y$ and each open $U \supseteq F$ there exists a $B \in \mathcal{B}$ with $F \subseteq B \subseteq U$. *I.e.*, a hyperbasis is a basis for the neighborhoods of the closed sets. Let the *hyperweight*, hX, of a space X be the least cardinality of an hyperbasis for X. As usual, wX and $\chi$X identify the weight and character of X.

In the next result, part (d) exposes a fact central to many considerations of the normality of finite products and part (e) can be viewed as the core (along with Stone's theorem that $N^\Omega$ is not normal) of most proofs demonstrating compactness properties of uncountable normal products.

**4.1 Proposition.** Let X be any space and Y the set $\kappa^+$ with a test worthy topology.
(a) If $\chi Y \leq \kappa$ and X is $[\omega,\kappa]$-compact, then $\pi_Y: X \times Y \to Y$ is closed.
(b) If $\pi_Y$ is closed, then X is $[\omega,\kappa]$-compact.
(c) If X and Y are normal with X $\kappa$-paracompact, $h(Y) \leq \kappa$, and $\pi_X: X \times Y \to X$ closed, then $X \times Y$ is normal.
(d) If $X \times Y$ is normal and $\pi_X: X \times Y \to X$ is z-closed, then X is $\kappa$-paracompact.
(e) If X is pseudocompact normal then $X \times \overline{\kappa}$ is normal iff X is $[\omega,\kappa]$-compact.

**Proof.** (a) The projection along an $[\omega,\kappa]$-compact space X into a space each point of which has a neighborhood basis of cardinality at most $\kappa$ is always closed: if $X \times \{y\}$ does not meet a closed set F it has a $\kappa$-fold cover $\{U_\alpha \times V_\alpha\}$ by basic open sets which do not meet F and hence a finite such cover, so y has a neighborhood V such that $X \times V$ does not meet F.

(b). Suppose $\pi_Y$ is closed but X is not $[\omega,\kappa]$-compact, that $\mathcal{O}$ is an open cover of X with no finite subcover, and that $\mathcal{O}$ is of minimum cardinality which (because $\overline{\lambda}$ and thus $\pi: X \times \overline{\lambda} \to \overline{\lambda}$ is closed for each cardinal $\lambda < \kappa$) we may suppose is $\kappa$. Then for F the "staircase" described in 1.3, the projection, $[0,\kappa)$, of F onto $[0,\kappa]$ is not closed since $\kappa$ is not isolated. Therefore X is $[\omega,\kappa]$-compact.

(c). For F and H disjoint closed subsets of X x Y and x∈X set $F_x = \pi_Y(F \cap (\{x\} \times Y))$ and $H_x = \pi_Y(H \cap (\{x\} \times Y))$. More generally, set $F_S = F \cap (S \times Y)$ and $H_S = H \cap (S \times Y)$ for any $S \subseteq X$. Let $\mathcal{B}$ be an hyperbasis for Y with $|\mathcal{B}| \leq \kappa$ which, we may suppose, includes Y and the empty set. For each ordered pair (B,C) of members of $\mathcal{B}$ which have disjoint closures set $O_{BC} = \{x \in X : \exists$ a neighborhood U of x with $F_U \subseteq U \times B$ and $H_U \subseteq U \times C\}$.

Since it contains a neighborhood of each of its points each $O_{BC}$ is open; we show that $\{O_{BC}\}$ covers X. For x in X, if $F_x$ is empty x is in $X \setminus \pi_X F = O_{\varnothing X}$, and if $H_x$ is empty x is in $X \setminus \pi_X H = O_{X\varnothing}$, which are open since $\pi_X$ is closed. Otherwise, since Y is normal it includes open subsets $V_0$, $V_1$, $W_0$ and $W_1$ with $F_x \subseteq V_0 \subseteq \overline{V_0} \subseteq V_1$, $H_x \subseteq W_0 \subseteq \overline{W_0} \subseteq W_1$, and $\overline{V_1} \cap \overline{W_1}$ empty. Set $S = Y \setminus W_0$, $T = Y \setminus V_0$, $T_F = \pi_X(F \cap (X \times T))$, and $S_H = \pi_X(H \cap (X \times S))$.

Since S, T, and $\pi_X$ are closed, $U = (X \setminus S_F) \cap (X \setminus T_H)$ is open. Notice that $F_U \subseteq U \times V_0$ and $H_U \subseteq U \times W_0$, so for B and C in $\mathcal{B}$ with $\overline{V_0} \subseteq B \subseteq V_1$ and $\overline{W_0} \subseteq C \subseteq W_1$, $F_U \subseteq U \times B$ and $H_U \subseteq U \times C$. Furthermore, since $\overline{V_1} \cap \overline{W_1} = \varnothing$, B and C have disjoint closures as required by the definition of $O_{BC}$. Therefore x is in $O_{BC}$.

Since $|\mathcal{B}| = \kappa$ it is a $\kappa$-fold open cover which we re-index as $\{O_\alpha\}$, labeling the corresponding members of $\mathcal{B}$ as $B_\alpha$ and $C_\alpha$ (so, for example, $F \cap (O_\alpha \times Y) \subseteq O_\alpha \times B$). Since X is $\kappa$-paracompact normal $\{O_\alpha\}$ has, by 2.2.3, a locally finite shrinking $\{V_\alpha\}$. Since $\{V_\alpha\}$ is a cover of X, F is contained in $U = \bigcup \{V_\alpha \times B_\alpha\}$, and since any point in $\overline{U}$ must be in the union of finitely many of the sets $\overline{V_\alpha} \times \overline{B_\alpha}$ which are disjoint from H, $\overline{U} \cap H = \varnothing$. Thus U and $(X \times Y) \setminus \overline{U}$ separate F and H and consequently $X \times Y$ is normal.

(d). Suppose X is not $\kappa$-paracompact; by 2.2.2(d) there exists a well indexed $\kappa$-fold open cover $\mathcal{O}$ of of X which has no locally finite shrinking, Since Y is test worthy our hypotheses apply to cardinals smaller than $\kappa$, so we may assume that k is minimal, *i.e.*, that $\mathcal{O}$ contains no repetitions. For $F = \bigcup_\alpha (F_\alpha \times [0,\alpha])$ the staircase set of 1.3, $F_\alpha = X \setminus O_\alpha$, and $H = X \times \{\kappa\}$, F and H are disjoint closed sets so by normality there is a Urysohn function f in $C(X \times Y, \mathbf{I})$ with $f(F) = \{0\}$ and $f(H) = \{1\}$.

Set $U_{\alpha n} = f^{-1}((1/n, 1]) \cap (X \times \{\alpha\})$; then $U_{\alpha n} \times \{\alpha\} \subseteq O_\alpha \times \{\alpha\}$ since f is zero on $F_\alpha \times \{\alpha\}$. Also, the collection $\mathcal{V} = \{U_{\alpha n}\}$ is a cover since $f(x,\kappa) = 1$ and hence (because f is continuous and $\kappa$ is not isolated) $f(x,\alpha)$ must be greater than zero for some $\alpha$. Thus $\mathcal{V}$ refines $\mathcal{O}$ and is, in fact, a shrinking. By 2.2.6 $\mathcal{V}$ has an open locally finite refinement, and that refinement is also a shrinking. This contradicts our assumption and shows that X is $\kappa$-paracompact.

(e) Since $h(\overline{\kappa}) = \kappa$ and the projection along $\overline{\kappa}$ is closed, $[\omega,\kappa]$-compactness (indeed, $\kappa$-paracompactness) of X implies, by (c), that $X \times \overline{\kappa}$ is normal. On the other hand, since $X \times \overline{\kappa}$ is pseudocompact (as is any product of a pseudocompact space with a compact space), normality of $X \times \overline{\kappa}$ implies, by 2.2.4, that X is $[\omega,\kappa]$-compact.

Statement 4.1(b) is due, in essence, to S. Mrówka, [1959, Remark, p21]. In establishing his fundamental characterization of compactness: X is compact if and only if the projection $\pi_Y$: X x Y→Y is closed for each space Y, he showed, for the "only if" direction, that X is compact provided $\pi$: X x $\overline{\kappa}$ →$\overline{\kappa}$ is closed for large enough cardinals κ.

Statement 4.1(c) is a special case of Noble [1971, Theorem 2.4] which extended Morita [1962, Theorem 2.2], the case with Y compact and hY replaced by wY, its weight.

Statement 4.1(d), as it applies to $\overline{\kappa}$, is due, for the countable case, to Dowker [1951]; his proof that for X normal X x **I** is normal if and only if X is countably paracompact actually only uses, in the only if direction, that X x **N*** is normal. For the general case, Morita [1962] showed, as noted above, that X x **I**(κ) is normal for each Tychonoff cube **I**(κ) if and only if X is κ-paracompact, a result he soon (in the subsequently submitted but more quickly published Morita [1961]) extended to Cantor and other cubes of complete metric spaces. His [1962, Lemma 2.5] suggests the stronger (in one direction) 4.1(d) with Y = $\overline{\kappa}$, which result is attributed to K. Kunen ( basis not indicated ) by both Katuta [1977] and Przymusinski [1984].

**4.2 Corollary**. Suppose Y = $\Pi_{\alpha \in A} Y_\alpha$ where $|A| = \kappa$ is infinite. Then:
(a) if no $Y_\alpha$ is indiscrete and $\pi$: X x Y→Y is closed then X is [ω,κ]-compact;
(b) if no $Y_\alpha$ is vacuously normal and X x Y is normal then X is κ-paracompact; hence
(c) if $X^\kappa$ is normal with κ uncountable, X is vacuously normal or [ω,κ]-compact; and thus
(d) if $X^\kappa$ is normal for uncountable κ ≥ wX, X is either vacuously normal or compact.

**Proof.** (a). Since not indiscrete each $Y_\alpha$ contains a nonempty proper closed subset $S_\alpha$. Let $s_\alpha$ be a point in $S_\alpha$, $t_\alpha$ a point in $X \setminus S_\alpha$, and $T_\alpha = \{s_\alpha, t_\alpha\}$. Each $T_\alpha$ is homeomorphic to either **II** or **2**; suppose κ of them are homeomorphic to **II** and set Z = $\Pi_\alpha Z_\alpha$ where $Z_\alpha = T_\alpha$ if $T_\alpha \cong$ **II** and $Z_\alpha = \{s_\alpha\}$ otherwise. Note that since X x Z = $\pi^{-1}$Z, the restriction of the closed map $\pi$ to X x Z is closed. Since Z is homeomorphic to **A**(κ), it contains a subset homeomorphic to ↑$\overline{\kappa}$ and thus $\pi$: X x ↑$\overline{\kappa}$ →↑$\overline{\kappa}$ is closed, so by 4.1(b) X is [ω,κ]-compact. If, on the other hand, κ of the $Z_\alpha$ are isomorphic to **2** we can construct Z homeomorphic to **C**(κ) and conclude that $\pi$: X x $\overline{\kappa}$ →$\overline{\kappa}$ is closed and hence that X is [ω,κ]-compact.

(b). By 1.5.5 X x $\overline{\kappa}$ is normal (with $\pi_X$ closed) so by 4.1(d) X is κ-paracompact.

(c). If X is not loosely compact it contains an infinite discrete family of nonempty irreducible closed subsets $\{S_n\}$ with the quotient q: S = $\bigcup_n S_n$ →**N** defined by $q(S_n) = \{n\}$ continuous and very open, causing $q^\kappa$: $S^\kappa$ →$N^\kappa$ to be continuous and very open. Since S must be closed, it is normal, implying that $N^\kappa$ is normal and thus (by Stone's theorem) that κ is countable. Thus X is loosely compact normal; by (b) it is κ-paracompact, so by 2.2.4 it is [ω,κ]-compact.

(d). If not vacuously normal X is, by (c), [ω,κ]-compact. Since each open cover has a refinement by basic open sets which must be of cardinality less than or equal to wX, hence less than or equal to κ, it follows that X is compact.

It is convenient, when dealing with the conventional product topology, to utilize a test space which occurs as a subspace of such products and "tests" all relevant cardinals at once. An alternative, convenient for testing the normality of box products, is to use a <u>collection</u> of test spaces, in which case the "test worthy" condition can be relaxed or dropped. For such testing, spaces such as those we will see in 6.3 with exactly one non-isolated point are particularly useful. These spaces are always normal: of any two disjoint closed subsets, one is open; and there is a well developed understanding of when the projections on such spaces are z-closed. There is also an extensive literature related to paracompactness and normality of box products which will not be surveyed here.

## 5 Infinite normal products

In contrast to the behavior of finite products, normality of countable products is about as agreeable as possible. (Assuming, of course, that the finite sub-products behave.)

**5.1 Theorem.** Suppose $X = \Pi_{n \in \mathbf{N}} X_n$ where no $X_n$ is vacuously normal. If X is normal it is countably paracompact, and if X is countably paracompact with each finite sub-product normal, then X is normal.

Zenor [1971, Theorem A], announced in [1969], states the corresponding result for countable powers of an Hausdorff space. Nagami - announced in [1969], presented in [1972] - proves a corresponding result for inverse limits of Hausdorff spaces and both authors also treat, by similar means, properties other than normality. Zenor's proof, with only the obvious modifications, gives 5.1. His proof relies on a lemma which I will replace with a perhaps more familiar result, the core of Tychonoff's proof that regular Lindelöf spaces are normal (see, *e.g.*, Willard [1970, Theorem 16.8] ):

**5.2 Tychonoff's Lemma.** If $\{U_n\}$ and $\{V_n\}$ are sequences of open sets which cover subsets F and H respectively with, for each n, $\overline{U}_n \cap H = \varnothing = F \cap \overline{V}_n$, then F and H are contained in disjoint open sets.

In tandem with Tychonoff's construction we will employ a property of countable products which is a recurring theme in proofs about them: if U is an open subset of $X = \Pi_n X_n$, then U can be represented as $U = \bigcup_n ( U_n \times X^{>n} )$, the union of an increasing sequence of open product sets of finite span.  It is somewhat simpler to construct the dual: $F=\bigcap_n(F_n \times X^{>n})$ where $F_n = cl(\pi_{\leq n}F)$.  (Notation: $\pi_{\leq n}$ projects X onto $\Pi_{i \leq n}X_i = X^{\leq n}$; $X^{>n} = \Pi^{i>n}X_i$.)  This representation has been widely exploited in the context of open covers - see Noble [2019, Sections 7 & 8] for some history and examples.

**Proof of 5.1**.  First suppose $X = \Pi_n X_n$ is countably paracompact with each finite subproduct normal, let F and H be disjoint closed subsets of X, and (as above) represent them as $F = \bigcap_n(F_n \times X^{>n})$ and $H = \bigcap_n(H_n \times X^{>n})$ where $F_n = cl(\pi_{\leq n}F)$ and $H_n = cl(\pi_{\leq n}H)$).  Set $G_n = F_n \cap H_n$; we want to separate F and H; if any $G_n$ is empty there exists open sets U and V which separate $F_n$ and $H_n$ ( since $X^{\leq n}$ is normal ) and $U \times X^{>n}$ and $V \times X^{>n}$ then separate F and H.  Thus we may suppose each $G_n$ is nonempty.  Since $\bigcap_n(G_n \times X^{>n})$ is empty ( any point in that intersection would be in $F \cap H$ ) countable paracompactness implies that this decreasing family of closed sets has an open expansion $\{E_n\}$ with $\bigcap_n \overline{E}_n$ empty where for each n, $G_n \times X^{>n} \subseteq E_n$.  Let $W_n = \pi_{\leq n}E_n$: since each point has a product set neighborhood which misses an $\overline{E}_n$ and thus is not in the corresponding $W_n \times X^{>n}$, the intersection $\bigcap_n(W_n \times X^{>n})$ is empty,  Also, $F_n \setminus W_n$ and $H_n \cup \overline{W}_{n+1}$ are disjoint closed subsets of the normal space $X^{\leq n}$ so there exists open sets $U_n \supseteq F_n \setminus W_n$ with $\overline{U}_n \cap (H_n \cup \overline{W}_{n+1}) = \emptyset$.

Since each point of F is in each $F_n \times X^{>n}$ but is not in some $W_n \times X^{>n}$, $F \subseteq \bigcup_n U_n$, and since each point of H is in each $(H_n \cup \overline{W}_{n+1}) \times X^{>n}$, $\overline{U}_n \cap H = \emptyset$.  Similarly, H can be covered by a family of open subsets $\{V_n\}$ whose closures do not meet F, so by Tychonoff's Lemma F and H can be separated by disjoint open sets, and therefore X is normal.

Next suppose X is normal and $\{G_n\}$ is a decreasing sequence of closed subsets with empty intersection; we will show that X is countably paracompact by constructing an expansion for $\{G_n\}$.  The factors of X are not vacuously normal, so for each n there exist points $p_n$ and $q_n$ with $\overline{p}_n$ and $\overline{q}_n$ disjoint: set $P^n = \Pi_{i>n}\overline{p}_n$, $Q^n = \Pi_{i>n}\overline{q}_n$, $F_n = cl(\pi_{\leq n}G_n)$, $F = \bigcup_n(F_n \times P^n)$, and $H = \bigcup_n(F_n \times Q^n)$.  Since $P^n$ and $Q^n$ are disjoint, F and H are disjoint.  Furthermore, F is closed: since $\bigcap_n P_n = \emptyset$, for any point z of X there exist an integer j such that z is not in $P_n$ for $n \geq j$ and an integer k, which we may assume is greater than or equal to j, such that a neighborhood $S \times X^{>k}$ of z does not meet $P_n$ for $n \geq k$.  That neighborhood does not meet $F_n \times X^{>n}$ for $n \geq k$, so z is in the closure of F only if it is in the closure of $\bigcup_{n<k}F_n \times P^n$, *i.e.*, only if it is in F.  Similarly, H is closed and hence, by the normality of X, there exist open sets $U \supseteq F$ and $V \supseteq H$ with $\overline{U} \cap \overline{V} = \emptyset$.  Note that, for $W_n = (\pi_{\leq n}U \cap \pi_{\leq n}V) \times X^{>n}$, $F_n \subseteq \pi_{\leq n}U \cap \pi_{\leq n}V$ so $G_n \subseteq W_n$ and $\bigcap_n \overline{W}_n \subseteq \overline{U} \cap \overline{V} = \emptyset$.  Thus $\{W_n\}$ form an expansion of $\{G_n\}$ and therefore X is countably paracompact.

The **Lindelöf number** $L(X)$ of a space X is the smallest infinite cardinal $\lambda$ such that each open cover of X has a subcover of cardinality at most $\lambda$; thus if $L(X) = \omega$ the space X is Lindelöf and possibly compact. Since each open cover of a space X has a subcover of cardinality no greater than the weight of X (one identified by a refinement composed of members of a base of that cardinality), $L(X) \leq w(X)$.

**5.3 Theorem.** If $X = \Pi_{\alpha \in A} X_\alpha$ is normal with $\kappa$ of its factors not vacuously normal and $\lambda$ is the cardinal such that, with at most countably many exceptions, each factor is either vacuously normal or has Lindelöf number less than or equal to $\lambda$, the index set A can be expressed as the union of disjoint subsets B, C, D, and E (some possibly empty) where:
- $X_B$ is vacuously normal;
- $X_C$ is compact;
- $X_D$ is $[\omega,\kappa]$-compact and D is empty if $\kappa \leq \omega$ or $\kappa \geq \lambda$; and
- $X_E$ is $\kappa$-paracompact with E countable.

**Proof.** Set $B = \{\alpha \in A: X_\alpha$ is vacuously normal $\}$ and $C = \{\alpha \in A \setminus B: X_\alpha$ is compact $\}$; $X_B$ is vacuously normal by 3.1(c) and $X_C$ is of course compact. If $\kappa$ is countable we will take D to be empty, so consider the case with $\kappa$ uncountable. By 2.1.5 all but at most countably many of the factors indexed by $A \setminus B$ are loosely compact and by 1.5.5 each such factor $X_\alpha$ has $X_\alpha \times \overline{\kappa}$ normal. It follows by 2.1.4, since each such $X_\alpha \times \overline{\kappa}$ is also pseudocompact, that the projections $X_\alpha \times \overline{\kappa}$ are closed and hence by 4.2(a) that each such factor is $[\omega,\kappa]$-compact. It follows, Noble, [1967; Theorem 4.1], that some product of all but countably many of the factors indexed by $A \setminus B$ is $[\omega,\kappa]$-compact, so there exists a subset $D' \subseteq A \setminus (B \cup C)$ with $X_D$ $[\omega,\kappa]$-compact and $E' = A \setminus (B \cup C \cup D')$ countable. If $\kappa \geq \lambda$ each $[\omega,\kappa]$-compact factor with a Lindelöf number less than or equal to $\lambda$ is compact so D' is countable and we set $D = \emptyset$ and $E = D' \cup E'$; otherwise $D' = D$ and $E' = E$ are as desired. Finally, $X_E$ is $\kappa$-paracompact: trivially if $\kappa$ is finite, otherwise by 4.2(b).

By citing 2.1.4 and 4.2(a) this proof uses, in essence, Mrówka's theorem that projections closed implies compactness. As an alternative one could use 4.2(b) plus 2.2.4 to show *κ*-paracompactness hence $[\omega,\kappa]$-compactness, in essence using Morita's theorem that normality of products with Tychonoff cubes implies paracompactness. Either approach leads to the desired conclusion that the factors are $[\omega,\kappa]$-compact. Such alternative proofs can, of course, shed light on different aspects of a topic. Indeed, having resumed the study of point set topology after a prolonged absence, I often find myself focusing more on the "commonality" of proofs than on the sometimes disconcerting plethora of nuanced properties which have been developed. (To which plethora I continue to contribute.) Thus motivated, I examine below several additional alternatives to the proof that a $T_1$ space all of whose powers are normal must be compact.

# 6  Alternative proofs that  X $T_1$ and $X^\kappa$ normal $\forall \kappa$ implies X compact

The original proof - Noble [1971, Corollary 2.2] - was based upon the implication of Stone's theorem that X must be countably compact plus an assortment of other results. It prompted a series of improvements, all also starting with Stone's theorem. Keesling [1972] gave two proofs and others followed: Franklin & Walker [1972], Polkowski [1979], and Przymusinskii [1984, Corollary 6.6]. As with the original proof, each of these relied upon some substantial additional theorem from the literature. Finally, Engelking [1988] gave a direct proof (still relying upon Stone's theorem) exhibiting two closed subsets in a power of X which cannot, if X is countably compact, be separated by open sets. Although I have, in a loose sense, kept my original approach, I have drawn from these proofs in crafting the presentation above. In particular, use of the Cantor cube formalism and the theorems of Morita and Mrówka was suggested by Keesling's second proof and the proofs of Polkowski and Przymusinski.

Before examining these alternative approaches I want to consider some other proofs which form the foundation upon which they stand:

- some alternate proofs of Stone's theorem that $\mathbf{N}^\Omega$ is not normal;
- Bourbaki's proof that the projection along a compact space is closed; and
- Bourbaki's proof that a space along which any projection is closed is compact.

As mentioned above, Stone's theorem is used in each of the alternates we will discuss to the proof that an uncountable product of non-compact $T_1$ spaces is not normal. The second Bourbaki proof is significant here because it underlies many (perhaps most) of the other proofs showing that particular product spaces are not normal.

## 6.1 Alternate proofs of Stone's theorem

Stone's proof considers disjoint closed subsets $F_0$ and $F_1$, the set of functions f: $\Omega \to \mathbf{N}$ which are one-to-one except on $f^{-1}(0)$, and the corresponding set of functions one-to-one except on $f^{-1}(1)$. As noted above, Stone proved directly that $F_0$ and $F_1$ cannot be separated. Referencing a theorem by Bockstein [1948] that disjoint open subsets of a product of separable metric spaces remain disjoint when projected into a suitably chosen countable sub-product, Corson [1959] gave a concise indirect proof of the fact that $F_0$ and $F_1$ cannot be separated (and thus of Stone's theorem) by noting that their projections onto any sub-product indexed by $\{\alpha_n\}$ will both contain the function p with $p(\alpha_n)=n$.

Bockstein's theorem is part of a body of results which include theorems guaranteeing that each continuous real-valued function on a particular product will depend on countably many coordinates, that is, will equal $g \circ \pi$ for some continuous g where $\pi$ is a projection into a countable sub-product. Comfort & Negropontis [1982, Chapter 10] contains a summary, extension, and history of such results. Mentioning Corson's proof, and noting that in a product of separable spaces any continuous function with regular second countable range depends on countably many coordinates, Ross & Stone [1964] give a third proof that $F_0$ and $F_1$ cannot be separated: if $\pi$ projects onto the sub-product indexed by $\{\alpha_n\}$ then no $g \circ \pi$ can be a Urysohn function separating $F_0$ and $F_1$ since for $p(\alpha_n)=n$ as above, $g^{-1}(p)$ is in both $F_0$ and $F_1$.

Each of these proofs of Stone's theorem use the same two closed subsets, $F_0$ and $F_1$. Fox [2013] explores the question as to what other closed subsets of $\mathbf{N}^\Omega$ might be used.. Among other results, Fox demonstrates that:

- There exist two countable closed discrete subsets of $\mathbf{N}^\Omega$ which cannot be separated by open sets having disjoint closures (as must be possible in a normal space); and
- given any nonempty closed subset F of $\mathbf{N}^\Omega$ with non-Lindelöf boundary, there exists a countable closed discrete subspace H which cannot be separated from F.

Note that $F_0$ and $F_1$ have empty interior and are not Lindelöf ($\{\pi_\alpha^{-1}(0)\}$ is an open cover of $F_0$ with no countable subcover), and recall that, in a normal space, not only can disjoint nonempty closed subsets be separated by open subsets, they can be separated by open subsets which have disjoint closures.

**6.2 Bourbaki** [1940, Ch.I, Sec. 10, Ex. 8, p68];  X compact implies $\pi_Y$ closed for all Y.

Engelking [1989, p133] describes the history of this result as follows: "Kuratowski proved in [1931] that in the class of metric spaces projections parallel to a compact space are closed mappings; this was generalized to topological spaces by Bourbaki in [1940] ....". Once a definition of compactness was settled the proof was straightforward: for p in the boundary of $\pi_Y F$, cover X x {p} with few enough open sets $\{U_\alpha \times V_\alpha\}$ that $V = \bigcap_\alpha V_\alpha$ is open and note that V is a neighborhood of p which does not meet $\pi_Y F$. It easily generalizes in two ways: For X κ-compact (aka [κ,∞]-compact, each open cover has a subcover of cardinality less than κ) and Y <κ-discrete (each intersection of fewer than κ open sets open) and for X [ω,κ]-compact and Y a space with a κ-fold base.

Those results suggest that, more generally, $\pi_Y$: X x Y → Y is closed if a class of instruments of convergence which in X must always have adherence is sufficient, in Y, to identify the closure of each subset. And that suggests consideration of a converse, our next topic.

**6.3 Bourbaki** [1961, Ch. I, Sec. 10.2, Lemma 1]: $\pi_Y$ closed for all Y implies X compact.

Bourbaki's proof (which does not appear in [1940] or [1951]) proceeds, given a space X which is not compact, by extending it to a space $X^* = X \cup \{\mathcal{F}\}$ such that the closure of the diagonal of X x X* projects to X, which is not closed. The construction: for $\mathcal{F}$ a filter base of closed subsets of X and for $S \subseteq X$ set $S^* = S \cup \{\mathcal{F}\}$. In X* take X to be discrete and let $\{F^*: F \in \mathcal{F}\}$ form a subbase for the neighborhoods of $\mathcal{F}$. Claim: if the projection along X into X* is closed, then in X the intersection of $\mathcal{F}$ is not empty; hence X is compact.

**Proof.** Let $\Delta = \{\{x,x\}: x \in X\}$; if the projection of $\overline{\Delta}$ onto X* is closed, it must contain $\mathcal{F}$, *i.e.*, there must exist a point $x \in X$ such that $(x, \mathcal{F})$ is in $\overline{\Delta}$. It follows that x is in each member of $\mathcal{F}$ since if $x \notin F$ then $(X \setminus F) \times F^*$ is a neighborhood of $(x, \mathcal{F})$ which does not meet $\Delta$.

We will see this construction again in 6.5 and 6.10. Notice that for this argument to hold X as a subspace of X* need not be discrete and the sets F* need not be open; all that is required is that the members of $\mathcal{F}$ be closed in the relative topology on X and that the added point $\mathcal{F}$ not be isolated. In particular, if X is dense in a space Y and y is a point in $Y \setminus X$ which has a neighborhood base $\mathcal{F}$ such that $\bigcap \{\overline{F}: F \in \mathcal{F}\}$ does not meet X, then the projection along X into $X^* = X \cup \{y\}$ is not closed.

From an internal viewpoint, if $\mathcal{F}$ is a filter base of cardinality $\kappa$ of closed subsets with empty intersection and X is $\kappa$-paracompact normal, then by 2.2.3 $\mathcal{F}$ has an expansion $\mathcal{G}$, a filter base of open subsets whose closures have empty intersection. Taking $\{G^*: G \in \mathcal{G}\}$ as a base for the neighborhoods of the point $\mathcal{F}$ in $X^* = X \cup \{\mathcal{F}\}$ produces a topology for X* for which the projection along X is not closed.

If, in our construction, the projection is closed, then the closed subset $H = X \times \{\mathcal{F}\}$ can be separated from any closed subset disjoint from it. If the projection is not closed it does not follow that H cannot be separated from such a set, but it suggests a likely place to find such an example, and thus an example for which the product X x X* is not normal. (Hence the staircase constructions of 1.3.) Replacing the requirement that the projection be closed with a weaker condition which still insures separation increases the likelihood of success in this endeavor.

In the case of separating $\overline{\Delta}$ from H we can replace the closed projection condition with the requirement that X x X* be C*-embedded in the subspace $Y = X \times X^* \cup \{(\mathcal{F}, \mathcal{F})\}$ of X* x X*. Specifically, if X x X* is normal then there exists a Urysohn function which is zero on $\overline{\Delta}$ and one on H, and such a function cannot extend to be continuous at $(\mathcal{F}, \mathcal{F})$ since that point is in the closure of both sets. This observation has implications for finite products, but also for infinite self products: indexing $\mathcal{F}$ as $\{F_\alpha: \alpha \in \kappa\}$, consider $F = \Pi_\alpha F_\alpha$ and the diagonal of the product of $\kappa$ copies of X. This approach is the basis for Franklin and Walker's proof which will be discussed in 6.7.

**6.4 Noble (1971):** if all powers of a $T_1$ space X are normal, X is compact.

The original proof of this fact documented the process by which I (primed by a question posed by Tony Hager) discovered the result. It suffices to show that the projections along each subproduct are closed since then, by Theorem 1.2 ( of Noble [1971] ) X is $[\omega,\kappa]$-compact. By Stone's theorem X is countably compact, hence pseudocompact, so by a theorem of Tamaño, or alternatively by a result concerning the Exponential Law for function spaces in Noble [1969b][5], the projections on X are z-closed which, as observed by Hager [1969], implies for X normal that its projections are closed. Curiously, a few authors sought a more direct proof.

**6.5 Keesling's first proof (1972)**

As with the sets discussed in 6.3, Keesling starts with a filter base $\mathcal{F}$ of closed subsets of X with empty intersection, the diagonal $\Delta$ of $\Pi\{X_Z: Z\in\mathcal{F}\}$ with each $X_Z$ a copy of X, and the closed set $F=\Pi\{F_Z\}$ where each $F_Z$ equals Z. To show that they cannot be separated he proves, as a lemma, that a continuous real-valued function on a countably compact product space depends on countably many coordinates. Keesling states that the result is known but that he is not aware of an adequate reference, but I think he may have been too modest. I have not found evidence that this particular result had appeared in the literature; Engelking [1966] had shown that real-valued continuous functions depend on countably many coordinates of a product with each finite sub-product Lindelöf by a proof which in fact holds (as noted in Noble & Ulmer [1972, Proposition 2.1]) if in each finite sub-product each uncountable set has an accumulation point, in particular if the finite sub-products are Lindelöf or countably compact. That and many earlier results are generalized by the paper contemporaneous with Keesling's, Noble & Ulmer [1972, Theorem 2.3].

Using Stone's theorem to conclude countable compactness, Keesling applies his lemma to a Urysohn function f which separates F and $\Delta$ (closed since X is Hausdorff) to identify a countable $B\subseteq\mathcal{F}$ and a continuous g such that $g\circ\pi_B=f$. As a countable filter base of closed sets in a countably compact space, there must exists a point p in $\bigcap_{n\in B}F_n$. But p would be in both $\pi_B(F)$ and $\pi_B(\Delta)$ which are disjoint since they are separated by f. This contradiction shows that no such $\mathcal{F}$ can exist, hence that X is compact.

---

5 With suitable separation Tamaño's theorem, [1960, Theorem 1], is in fact equivalent to the mentioned exponential law, Glicksberg's theorem [1959] characterizing when $\beta(X\times Y)=\beta X\times\beta Y$, Isbell's theorem [1964] describing when $X\times Y$ is fine (has a topology compatible with its finest uniformity); some versions of Ascoli's theorem; and many other formulations: Comfort & Hager [1971] list eighteen.

### 6.6 Keesling's second proof (1972)

Noting Morita's result: X is κ-paracompact normal if and only if $X \times \mathbf{I}(\kappa)$ is normal, Keesling argued, more or less as was done above in 1.2.2, that $X \times \mathbf{C}(\kappa)$ must be normal, so (by additional argument) $X \times \mathbf{I}(\kappa)$ must be normal and consequently (there being no limit on κ) X must be paracompact. But each paracompact lightly compact space is compact and X, being Hausdorff, is by Stone's theorem, lightly compact, so it is compact.

### 6.7 Franklin & Walker (1972)

Franklin and Walker use the construction discussed in 6.3 with $\mathcal{F} = \rho$ an z-ultrafilter in $\beta X \setminus X$. Thus $X^* = X \cup \{\rho\}$, $F = \Pi_{\alpha \in \rho} F_\alpha \subset \Pi_{\alpha \in \rho} X_\alpha \subseteq \Pi_{\alpha \in \rho} \beta X_\alpha$ where each $F_\alpha$ equal to its index, a zero set. F and the diagonal Δ are disjoint closed sets, but they cannot be separated: if f were a Urysohn function on $\Pi_\alpha X_\alpha$ separating them, then f would extend to a continuous function f* on $\Pi_\alpha \beta X_\alpha$ since by Stone's theorem $\Pi_\alpha X_\alpha$ is pseudocompact and therefore, by Glicksberg's theorem, equal to $\beta(\Pi_\alpha X_\alpha)$. But the point $(p_\alpha)$ with each $p_\alpha = \rho$ of $\Pi_\alpha \beta X_\alpha$ is in the closure of both F and Δ, so this is not possible. Therefore $\beta X = X$ and X is compact.

### 6.8 Polkowski (1979)

Recall that a continuous surjective function (a "map") is proper if it is also closed with each fiber compact. As was noted above Bourbaki showed that all powers of a proper map are proper, hence closed. Conversely, if all powers of a map are closed it is proper - Noble [1971, Corollary 3.7]. Observing that if all powers of X are normal $T_1$ then (by Stone's and Tamaño's theorems) the projection π: $X \times X \to X$ and all its powers are closed, Polkowski concludes that powers of π are proper and hence that X is compact.

To prove "directly" that a non-constant map with Hausdorff domain whose powers are closed must be perfect, Polkowski, arguing from Mrówka's theorem as discussed in Subsection 1.2, observes that it suffices to show for each Cantor cube $\mathbf{C}(\kappa)$ and for each fiber F of the map that the projection π: $F \times \mathbf{C}(\kappa) \to \mathbf{C}(\kappa)$ is closed. He does so by chasing circuits on a pair of commutative diagrams. Although presented differently, his approach is demonstrated by the next proof.

**6.8.1 Proposition.** Let m: $X \to Z$ be continuous and surjective with Z not indiscrete. All powers of m are closed if and only m is proper.
**Proof.** Products of proper maps are closed: Bourbaki [1961, p119]. Suppose all powers of m are closed and, for convenience, that Z is disjoint from X. We show that the fibers of m are compact. For $F = m^{-1}(z)$ such a fiber and $Z^0 = Z \setminus \{z\}$ set $Y = F \cup Z^0$; let $p: X \to Y$ be the identity on F and m on $X \setminus F$; and let $q: Y \to Z$ be the identity on $Z^0$ and m on F.

Note that $q \circ p = m$. Give Y the quotient topology; for H a closed subset of Z the set $q^{-1}(H)$ is closed in Y since p is a quotient map and $p^{-1}(q^{-1}(H)) = m^{-1}(H)$ is closed in X. Thus q is continuous. Also, since p is continuous the set map $r = p^{-1}$ is closed and thus, as m is also closed, $q = r \circ m$ is closed. It follows, since all powers of m are closed and all powers of p are continuous, that $r^\kappa$, $m^\kappa$, and thus $q^\kappa = r^\kappa \circ m^\kappa$, are continuous closed maps.

Consider $q^2: Y \times Y \to Z^2$ and label as $\pi$ the restriction of $q^2$ to the subset $F \times Y$; since $\pi(F \times Y) = \{z\} \times Y$ and $\pi_Y \circ \pi$ is the identity, $\pi$ is the projection. Thus the projection $\pi^\kappa$ on $F^\kappa \times Y^\kappa \to Y^\kappa$ is closed. Since Z is not indiscrete Y is not indiscrete so 4.2(a) applies to show that F is, for any infinite $\kappa$, $[\omega,\kappa]$-compact. Therefore F is compact, as desired.

### 6.9 Przymusinskii (1984)

The class of perfect pre-images of metrizable spaces coincides with the class of "paracompact p-spaces"; it is closed (as is the class of metrizable spaces) under countable products. Considering only Hausdorff spaces, Przymusinskii uses this fact, Stone's theorem, and the fact that the product of a paracompact space with a compact space is paracompact (due to Dieudonné and a consequence of Morita's theorem) to show that a product of paracompact p-spaces is paracompact (or normal) if and only if all but countably many of the factors are compact. Noting that if the product of X with a product having $\kappa$ factors is normal its product with the Cantor cube $C(\kappa)$ must be normal, Przymusinskii then invokes Morita's theorem to conclude that if all powers of X are normal X must be compact. I suspect this proof extends to one involving the class of proper pre-images of pseudometrizable spaces, which also are closed under countable products, but have not explored that possibility.

### 6.10 Engelking (1988)

Finally, Engelking provided an elementary proof derived from first principals (and Stone's theorem), exhibiting two closed subsets in a power of X which cannot, if X is countably compact, be separated by open sets. As in Bourbaki's 1981 proof, Keesling's First proof, and the proof by Franklin & Walker, Engelking uses a filter base $\{F_\alpha\}$ of nonempty closed sets with empty intersection, and constructs $F = \Pi_\alpha F_\alpha$. He shows, without extending X, that his F cannot be separated from the diagonal which, since X is Hausdorff, is closed. Engelking does so by constructing a sequence of points in F for which any cluster point gives rise to a point on the diagonal each neighborhood of which must meet F. He then invokes Stone's theorem to note that X is countably compact so such a point must exist,

## 7  Implications beyond point set topology

I warn the reader that I am not qualified to comment on the topics discussed in this section.

### 7.1 Categorical implications

Even before its publication in 1971 the fact that all powers of a $T_1$ space being normal implies compactness was applied by Herrlich & Strecker [1971] - submitted 15 December 1970 - and by Franklin, Lutzer & Thomas [1971] to give categorical characterizations of **CH**, the category of compact Hausdorff spaces, as a subcategory of all Hausdorff spaces.  The  Herrlich & Strecker version identifies **CH** as the only full non-trivial subcategory which is both varietal (in the sense of Linton) and epireflective. Franklin, Lutzer & Thomas find it to be the only non-trivial, productive, left-fitting subcategory preserved by shrinks, in **TOP**, of $T_2$ extremal points. (Earlier Franklin & Thomas [1970] had noted that work by DeGroot translates to a description of **CH** as the only nontrivial, productive subcategory $\mathcal{C}$ of **TOP** which is preserved by closed epis and satisfies the condition: if $X \in Ob\mathcal{C}$ and $m^\#:F \to X$ is an extremal mono, then $F \in Ob\mathcal{C}$  if and only if $m^\#$ is a closed map.

Subsequently  Franklin, Lutzer & Thomas [1977] presented a unified theory of subcategory characterization and provided characterizations based only on global behavior of the (objects and) morphisms for the subcategories $T_0$, $T_1$, $T_2$, $T_3$, and $T_4$, (with $T_2 \subset T_3 \subset T_4$ ) as well as paracompact Hausdorff, locally compact Hausdorff, separable Hausdorff, separable metrizable, and first countable spaces.

Notably, characterization of regular spaces, normal spaces, and even compact spaces are absent, and I have not found them characterized in later work.  For example, Herrlich [1993] explores almost reflective subcategories, but notes that compact, compact $T_0$, and compact $T_1$ spaces form implicational subcategories of **TOP** which are neither almost reflective nor orthogonal.  Perhaps some of the results here will lead to different, interesting characterizations in the realm of poorly separated spaces.

### 7.2  Lattices, locals, frames, quantrales, topos, domains, and such

Interest in and research into poorly separated spaces has largely been motivated by the problem of modeling computation and providing a context for programing languages and their semantics.  Recursion and branching in the operation of computing algorithms are modeled using a time ("system clock") dependent ordering of the potential states of a virtual computing mechanism. With independent inputs at the bottom (to which, generally, a least element is attached representing the loading of the initial conditions) the order progresses upward to states which can be achieved as the algorithm is executed.

Transfinite recursion is discouraged because of its infinite cost and infinite delay, but the height of the state space is taken as countably infinite as an efficient means of treating very large finite numbers.

Initial elements, which we can represent as numbers, are necessarily of finite resolution, but it can be quite difficult to establish the resolution an algorithm will require - other than by running it first with inputs of much higher resolution. See for example the discussion in Blackmore & Peters [2007, in particular 3.2, 4.2, and 5.1]. Consequently, and for much the same reason that algebraic and topological completeness is desirable in other areas of mathematics, the inputs to algorithms are generally assumed to be real or complex numbers, or elements of some finite or countable product of these ordered fields. Thus the width of the state space is usually taken to be uncountable.

I don't really speak (or read) "lattice lingo" so can't claim confidence in my search, but I have not found any lattice treatments of separation which go beyond some interesting variations of regularity. Given its prominence in topology I expect that normality has some relevance in the computing context. After all, "closed" sets contain the states from which their elements can be produced, and "open" sets contain the states which can be produced from their elements. Thus disjoint closed sets (provided we ignore the attached start element) represent sets of states each not reached from a state in the other, and the existence of "open" sets containing them indicates that they will never join. Thus they can be viewed as representing sub-algorithms which could be run in parallel.

Rudin, M.E. & Starbird, M.,
 [1975] *The normality of products with a metric factor*, General Topology and Appl. **5**, 235-248.

Scarborough, C. T. & Stone, A.H.,
 [1966] *Products of nearly compact spaces*, Trans. Amer. Math. Soc. **124**, 131-147.

Stephenson Jr., R.M.,
 [1969] *Product spaces for which the Stone-Weierstrass theorem holds,* Proc Amer. Math. Soc, **21**, 284-288.
 [2004] *Pseudocompact Spaces*, Chapter d-7 in Hart, Nagata & Vaughan [2004], 177-181.

Stone, A. H.,
 [1948] *Paracompactness and product spaces*, Bull. Amer. Math. Soc. **54**, 977-982.
 [1960] *Hereditarily compact spaces*, Amer. J. Math. **82**, 900-916.

Steen, Lynn A. & Seebach, J. Arthur Jr.,
 [1970] *Counterexamples in Topology*, Holt, Rinehart and Winston.

Tamaño, H.,
 [1960] *A note on the pseudo-compactness of the product of two spaces*, Mem. Coll. Sei. Univ. Kyoto Ser. A Math. **33** (1960/61), 225-230.

Willard, Stephen,
 [1970] *General Topology*, Addison-Wesley.

Zenor, Phillip,
 [1969] *Countable paracompactness and normality in product spaces*, Proc. 1969 Houston Top. Conf. (1971), 53-55.
 [1971] *Countable paracompactness in product spaces*, Proc. Amer. Math. Soc. **30,** 199–201.
38